\documentclass{nseJournal}
\usepackage{verbatim}
\usepackage{float}
\usepackage{amsmath}
\usepackage{amssymb}
\usepackage{graphicx}
\usepackage{subcaption}

\usepackage{lineno}

\makeatletter
\def\blfootnote{\gdef\@thefnmark{}\@footnotetext}
\makeatother

\newcommand{\cred}[1]{{\color{black}{#1}}}

\begin{document}

\title{An Efficient Sweep-based Solver for the $S_{N}$ Equations on High-Order Meshes}

\addAuthor{\correspondingAuthor{T. S. Haut}}{a}
\addAuthor{P. G. Maginot}{b}
\addAuthor{V. Z. Tomov}{a}
\addAuthor{B. S. Southworth}{c}
\addAuthor{T. A. Brunner}{b}
\addAuthor{T. S. Bailey}{b}

\correspondingEmail{haut3@llnl.gov}

\addAffiliation{a}{Center for Applied Scientific Computing, Lawrence Livermore National Lab\\ P. O. Box 7000, Livermore, CA 94551}
\addAffiliation{b}{Design Physics Division, Lawrence Livermore National Lab\\ P. O. Box 7000, Livermore, CA 94551}
\addAffiliation{c}{Department of Applied Mathematics, University of Colorado\\ Engineering Center, ECOT 225, 526 UCB, Boulder, CO 80309-0526}

\addKeyword{transport}
\addKeyword{discrete ordinates}
\addKeyword{high-order}
\addKeyword{sweep}
\addKeyword{unstructured}

\author{
 \vspace{10mm} \\
  T. S. Haut,$^{*,a}$ 
  P. G. Maginot,$^b$
  V. Z. Tomov,$^a$ \\
   B. S. Southworth,$^c$
   T. A. Brunner,$^b$ and T. S. Bailey$^b$ 
   \\
   \\
   \normalsize{
    \begin{itshape}
    \begin{tabular}[t]{c}
       $^a$ Center for Applied Scientific Computing, Lawrence Livermore National Lab\\  P. O. Box 7000, Livermore, CA 94551  \\  \\
       $^b$ Design Physics Division, Lawrence Livermore National Lab\\ P. O. Box 7000, Livermore, CA 94551 \\  \\
       $^c$ Department of Applied Mathematics, University of Colorado\\ Engineering Center, ECOT 225, 526 UCB, Boulder, CO 80309-0526
    \end{tabular}
    \end{itshape}
    }
     \vspace{10mm}
    \\
    $^*$Email: \url{mailto:haut3@llnl.gov}
    }

 \date{
        \vspace{15mm}
        \begin{tabular}[t]{lc}
        Number of pages:& \pageref*{LastPage} \\
        Number of tables:& \totaltables \\
        Number of figures:& \totalfigures
        \end{tabular}
    }
 \clearpage\maketitle
    \thispagestyle{empty}	    
 
\begin{abstract}
We propose a graph-based sweep algorithm for solving the 
steady state, mono-energetic discrete ordinates on meshes
of high-order curved mesh elements. Our spatial discretization
consists of arbitrarily high-order discontinuous Galerkin finite
elements using upwinding at mesh element faces. To determine mesh
element sweep ordering, we define a directed, weighted graph whose
vertices correspond to mesh elements, and whose edges correspond to
mesh element upwind dependencies. This graph is made acyclic by removing
select edges in a way that approximately minimizes the sum of removed
edge weights. Once the set of removed edges is determined, transport
sweeps are performed by lagging the upwind dependency associated with
the removed edges. 
The proposed algorithm is tested on several 2D and 3D meshes composed of high-order curved mesh elements.
\end{abstract}

\section{Introduction}

Hydrodynamics discretizations that leverage high-order (HO) curved
meshes have become increasingly popular in recent years, due to several
mathematical and computational advantages they afford. For example,
HO meshes can more accurately model curved geometries, and consequently
simulations that use HO meshes can use significantly fewer mesh elements
to achieve an equivalent accuracy in comparison to low-order (i.e.
straight-edged) meshes. HO meshes are essential for achieving optimal
finite element convergence rates (for smooth problems) on domains
with curved boundaries or interfaces, in both the Lagrangian \cite{Dobrev2012,Dobrev2013,Tomov2016}
and arbitrary Lagrangian-Eulerian (ALE) \cite{Boscheri2016,Dobrev2018}
contexts. Additionally, HO meshes permit better spherical symmetry
preservation for implosion calculations in axisymmetric geometries
\cite{Dobrev2013}. 

In the context of radiation transport, HO methods have the potential
for reducing the prohibitively large number of spatial unknowns required
to resolve the radiation's seven-dimensional phase space (three in space,
two in angle, energy, and time). The need to reduce the number of thermal
radiation transport (TRT) unknowns is further magnified with the rise
of memory motion/communication-bound next-generation architectures
such as the graphics processing unit based machine Sierra \cite{supinski2017}.
The ability to solve the TRT equations on HO curved meshes is further
motivated when one is interested in the simulation of multi-physics
phenomena, such as inertial confinement fusion, using HO hydrodynamics
methods, such as those in \cite{Dobrev2012,Dobrev2013,Tomov2016,Boscheri2016,Dobrev2018}.
Such radiation-hydrodynamics simulations must preserve the equilibrium
diffusion limit, which requires the solution of the $S_{N}$ equations
on a mesh compatible with the HO hydrodynamics mesh.

One possibility to achieve the required radiation hydrodynamics coupling
would be to interpolate fields such as the radiation energy density
between a HO representation on the hydrodynamics mesh and a field
discretized with linear elements on a refined, low-order (LO) TRT
specific mesh. However, this would necessitate a significant increase
in TRT degrees-of-freedom to ``straighten out'' the HO hydrodynamics
meshes. Additionally, (i) no analysis currently exists proving that
such a scheme preserves the equilibrium diffusion limit, and (ii)
interpolating between LO and HO representations of the same field
will likely lead to decreased robustness, stability, and physics fidelity.

An alternative approach is to solve the TRT equations directly
on the HO mesh using a spatial discretization for the TRT equations
that is consistent with the thermodynamic variables of the hydrodynamics
discretization. In \cite{Dobrev2012,Dobrev2013,Dobrev2018}, the hydrodynamics'
thermodynamic variables are discretized in space by a discontinuous
Galerkin (DG) formulation. Solving the TRT equations directly on the
hydrodynamics mesh is most consistent with historical precedent for
robust radiation hydrodynamics coupling stability \cite{kull}. Further,
when the TRT unknowns and the hydrodynamics' thermodynamic unknowns
have identical centerings, it is likely that even better multi-physics
robustness, numerical stability, and accuracy can be achieved. In
particular, Woods and Palmer have shown that an upwinded DG spatial
discretization of the mono-energetic, steady-state, $S_{N}$ angular
discretization of the linear Boltzmann equation (i) yields HO convergence
\cite{woodsPalmer1,woodsPalmer2} and (ii) maintains the thick diffusion
limit \cite{Woods2017a,woodsPalmer2}. Spatial discretizations of
the TRT equations that preserve the thick diffusion limit are tantamount
to spatial discretizations that will preserve the equilibrium diffusion
limit when applied to time-dependent, multi-frequency, $S_{N}$ angular
approximations to the non-linear thermal radiative transfer equations
\cite{morel_radtran,adams_nowak}. Hereafter, we refer to the spatially
analytic mono-energetic, steady-state, $S_{N}$ angular discretizations
of the linear Boltzmann equation simply as the linear transport equations.

When source iteration \cite{lewis_book} is used to iteratively
solve for particle scattering, each iteration requires solving the linear
transport equation (independently) for each angle. Upwinded spatial
discretizations of the linear transport equations on product meshes of
orthogonal mesh elements (or on general 2d, convex meshes),
lead to a block lower-triangular set of equations, in some ordering,
for each angle. This system can be solved on an element-by-element
basis by doing a block forward solve of the corresponding matrix; in 
the the transport literature, this is referred to as a ``transport sweep.''
On meshes of non-orthogonal linear (straight-edged) mesh elements,
upwinded DG discretizations of the linear transport
equations do not necessarily form lower-triangular systems
of equations. 
\cred{
Systems with upper-triangular components are often produced by discretizations
on 3D linear unstructured grids \cite{sn_on_tets}, or by discretizations on
2D linear meshes that contain concave elements. }
Nevertheless, it is usually
possible to choose a mesh ordering that results in a nearly block
lower-triangular matrix whose inverse is approximated well in a transport
sweep \cite{McLendon2001,Pautz2002,Plimpton2005,Pautz2017}. Such
an algorithm is equivalent to an ordered block Gauss-Seidel iteration
on the upwinded DG equations. 

On HO curved meshes, an upwind discretization of the linear
transport equations can result in significantly more upper triangular
entries, compared to similar problems composed on linear mesh elements.
As an example, it is possible to have two neighboring HO mesh elements
that are both upwind of each other, even in 2D. Thus, in HO meshes,
strongly connected components (SCCs) of mesh elements are expected
to be more prevalent and larger than SCCs appearing in linear meshes.
Each SCC corresponds to a non-triangular component of the
spatially discretized linear transport equations. As mesh element
order increases, such SCCs likely grow in size and number, causing
convergence of a naive Gauss-Seidel iteration to suffer. To that end,
methods for determining a mesh element ordering that minimizes the
size of these SCCs takes on greater importance with HO curved meshes.

The main purpose of this paper is to develop and present a transport
sweep algorithm based on a quasi-optimal mesh-element ordering for
a block Gauss-Seidel iteration. In particular, we extend the graph-based
methods of \cite{McLendon2001,Pautz2002,Plimpton2005,Pautz2017} to
HO curved meshes in a manner that directly connects the graph to both
the mesh topology and HO DG spatial discretization. After defining
such a graph, a transport-sweep mesh-element ordering is determined
so that the ``optimal'' approximation of the inverse of the upwinded
DG discretization can be applied on an element-by-element basis. In
graph-theoretic terms, we construct the maximum acyclic subgraph by
removing edges in the dependency graph until the modified graph is
acyclic using a heuristic that (approximately) minimizes the weight
of removed edges. Performance of this graph based, approximate transport
sweep algorithm is then demonstrated on a variety of test problems, including
an annulus embedded in a square domain, reminiscent of the canonical
light-water reactor pin cell \cite{bell_glasstone}; a HO 2D-mesh derived
from a HO Lagrangian simulation of the ``triple point'' problem
\cite{Dobrev2012,triplePoint}; and a 3D extension of the 2D, single-mode,
Rayleigh-Taylor
instability problem presented in \cite{Dobrev2012}. In this series
of test problems, the iteration count of our Gauss-Seidel, quasi-optimal
element-by-element transport sweep solver is compared against the
number of iterations required for solution of the exact upwinded DG
equations via an approximate ideal restriction algebraic multigrid
solver (AIR) \cite{air1,air2,hypre2_14}. Finally, we compare the unknown
count, iteration count, and accuracy of our Gauss-Seidel iteration on the
HO mesh to refined meshes of linear mesh elements using a linear
upwinded DG transport discretization.

\section{Spatial Discretization of the Linear Transport Equations}

In this work, we focus on the solution of the steady-state, mono-energetic
discrete ordinates linear Boltzmann equations,
\begin{align}
\boldsymbol{\Omega}_{d}\cdot\nabla_{\boldsymbol{x}}\psi_{d}\left(\boldsymbol{x}\right)+\sigma_{t}\left(\boldsymbol{x}\right)\psi_{d}\left(\boldsymbol{x}\right) & =\frac{\sigma_{s}\left(\boldsymbol{x}\right)}{4\pi}\varphi\left(\boldsymbol{x}\right)+q_{d}\left(\boldsymbol{x}\right),\,\,\boldsymbol{x}\in\mathcal{\mathcal{D}}\nonumber \\
\psi_{d}\left(\boldsymbol{x}\right) & =\psi_{d,\text{inc}}\left(\boldsymbol{x}\right),\,\,\,\,\mathbf{x\in\partial}\mathcal{D}\,\,\,\,\text{and}\,\,\,\,\mathbf{n}\left(\boldsymbol{x}\right)\cdot\boldsymbol{\Omega}_{d}<0.\label{eq:linearSnTransport}
\end{align}
In (\ref{eq:linearSnTransport}), $\mathcal{D}$ is a domain in $dim$
dimensions and $\mathbf{\partial}\mathcal{D}$ is its boundary; $\psi_{d}\left(\boldsymbol{x}\right)$
is the angular flux $\left[\frac{n}{cm^{2}-sec-ster}\right]$ in the
direction of discrete ordinate $\boldsymbol{\Omega}_{d}$, with $\boldsymbol{\Omega}_{d}\in\mathbb{S}^{2}$;
the total macroscopic interaction cross section, $\sigma_{t}\left(\boldsymbol{x}\right)\left[cm^{-1}\right]$,
is the sum of the macroscopic absorption cross section, $\sigma_{a}\left(\boldsymbol{x}\right)\left[cm^{-1}\right]$,
and the scattering macroscopic cross section, $\sigma_{s}\left(\boldsymbol{x}\right)\left[cm^{-1}\right]$;
scattering is assumed to be isotropic; we have assumed the validitiy
of the discrete ordinates approximation to the scalar flux, $\varphi\left[\frac{{n}}{cm^{2}-sec}\right]$,
\begin{equation}
\varphi\left(\boldsymbol{x}\right)\approx\sum_{d=1}^{N_{\Omega}}\omega_{d}\psi_{d}\left(\boldsymbol{x}\right)\,\,;\label{eq:snScalarFlux}
\end{equation}
$q_{d}\left(\boldsymbol{x}\right)\,\left[\frac{{n}}{cm^{3}-sec-ster}\right]$
is a volumetric source in the direction of $\boldsymbol{\Omega}_{d}$;
$\psi_{d,inc}\left(\boldsymbol{x}\right)$ is a known, incident angular
flux on the boundary $\partial\mathcal{{D}}$ of domain $\mathcal{D}$;
and $\mathbf{n}\left(\boldsymbol{x}\right)$ is the outward directed
unit normal on $\partial\mathcal{{D}}$ at $\boldsymbol{{x}}$. Additionally,
we have chosen the set of angular quadrature weights and directions,
$\left\{ (\omega_{d},\mathbf{\Omega}_{d})\right\} _{1}^{N_{\Omega}}$,
such that all $\omega_{d}>0$ and $\sum_{d=1}^{N_{\Omega}}\omega_{d}=4\pi$. 

\paragraph*{Discrete representation of the HO mesh.}

Equation (\ref{eq:linearSnTransport}) is discretized over a HO curved
mesh. Let $\cup\kappa_{e}$ denote a tesselation of physical space
using HO mesh elements, $\{\kappa_{e}\}_{1}^{N_{\kappa}}$, each being
of order $r$. Consider a set of scalar Lagrangian finite element
basis functions $\{\hat{v}_{m}(\xi)\}_{1}^{N_{v}}$ defined on the
reference element $\hat{{\kappa}}$. For simplicial elements (triangles
and tetrahedrons), the basis $\{\hat{v}_{m}\}$ spans $P_{r}$, the
space of all polynomials of total degree at most $r$; for tensor
product elements (quadrilaterals and hexahedrals) $\{\hat{v}_{m}\}$
is chosen to be a basis for $Q_{r}$, the space of all polynomials
of degree at most $r$ in each variable. The shape of any element
$\kappa_{e}$ in the mesh is then fully described by a matrix $\mathbf{{x}_{e}}$
of size $dim\times N_{v}$ whose columns represent the coordinates
of the element control points (aka element degrees of freedom). Given
$\mathbf{{x}_{e}}$, the position of $\kappa_{e}$ (corresponding
to a point $\xi$ in $\hat{{\kappa}}$) is derived via the mapping
$\mathbf{T}_{e}:\hat{\kappa}\rightarrow\kappa_{e}$:

\begin{equation}
\mathbf{x_{e}(\xi)}=\mathbf{T}_{e}\left(\boldsymbol{\xi}\right)=\sum_{m=1}^{N_{v}}\mathbf{x}_{e,m}\hat{v}_{m}\left(\boldsymbol{\xi}\right),\,\,\,\,\mathbf{x}_{e,m}=\mathbf{T}_{e}\left(\boldsymbol{\xi}_{m}\right),\label{eq:elementTransformation}
\end{equation}
where we use $\mathbf{x}_{e,m}$ to denote the $m$-th column of $\mathbf{{x}_{e}}$.
When two elements $\kappa_{e}$ and $\kappa_{e'}$ share a common
mesh entity (vertex, edge, or face) then their control-point matrices
$\mathbf{{x}_{e}}$ and $\mathbf{{x}_{e'}}$ are not independent because,
to ensure mesh continuity, the descriptions of the common entity (through
$\mathbf{T}_{e}$ on one hand and $\mathbf{T}_{e'}$ on the other)
must coincide. This type of interdependence among mesh elements is
typically expressed by defining a global vector $\mathbf{{x}}$ of
control points (degrees of freedom) and a set of linear operators
$\mathbf{\mathbf{{x}\rightarrow}{x}_{e}}$ , one per mesh element,
that define/extract the local coordinates $\mathbf{{x}_{e}}$ from
the global vector $\mathbf{{x}}$. Thus, the global continuity of
the mesh is ensured for any value of $\mathbf{{x}}$.

\paragraph*{Discontinuous Galerkin formulation.}

We spatially discretize (\ref{eq:linearSnTransport}) using HO DG
approximations of $\psi_{d}\left(\mathbf{x}\right)$ and $\varphi\left(\mathbf{x}\right)$.
Let $\{\hat{u}_{m}(\xi)\}_{1}^{N_{u}}$ be a set of scalar Lagrangian
finite element basis functions defined on the reference element $\hat{{\kappa}}$.
As before, this basis spans $P_{s}$ or $Q_{s}$ for simplicial and
tensor elements, respectively, where $s$ is the DG finite element
order. For a given element $\kappa_{e}$ in physical space, the corresponding
basis functions, $\{u_{m}^{e}(\mathbf{{x}})\}_{1}^{N_{u}}$, are defined
through 
\begin{equation}
u_{m}^{e}\left(\mathbf{x_{e}(\boldsymbol{{\xi}})}\right)\equiv\hat{u}_{m}\left(\boldsymbol{\xi}\right),\,\,\,\,\boldsymbol{\xi}\in\hat{\kappa}.\label{eq:localBasis}
\end{equation}
Regardless of the particular mesh element, each $u_{m}^{e}$ satisfies
$u_{m}^{e}(\mathbf{x}_{n})=\delta_{m,n}$, where $\mathbf{{x}}_{n}$
is the nodal position of the $n$-th basis function. These basis functions
are used to discretize the angular flux and scalar flux on any element
$\kappa_{e}$ as 
\[
\widetilde{\psi}_{e}^{d}\left(\mathbf{x}\right)=\sum_{m=1}^{N_{u}}\psi_{e,m}^{d}u_{m}^{e}\left(\mathbf{x}\right),
\]
and 
\[
\widetilde{\varphi}_{e}\left(\mathbf{x}\right)=\sum_{m=1}^{N_{u}}\varphi_{e,m}u_{m}^{e}\left(\mathbf{x}\right).
\]
On a face $\Gamma_{e,e'}$ between neigboring mesh elements $e$ and
$e'$, upwinding is used to define the angular flux,

\[
\widetilde{\psi}_{e,e'}^{d}(\mathbf{x})=\begin{cases}
\widetilde{\psi}_{e}^{d}(\mathbf{x})\mathbf{\quad\text{{if} \,\ensuremath{\boldsymbol{\Omega}_{d}}}\cdot\mathbf{n_{\textrm{e}}}}(\mathbf{x})\geq0,\\
\widetilde{\psi}_{e'}^{d}(\mathbf{x})\mathbf{\quad\text{{if} \,\ensuremath{\boldsymbol{\Omega}_{d}}}\cdot\mathbf{n_{\textrm{e}}}}(\mathbf{x})<0,
\end{cases}
\]
where $\mathbf{\mathbf{{n}}}_{e}(\mathbf{{x}})$ is the outward directed
normal of $\kappa_{e}$ that point into $\kappa_{e'}.$ Following
the standard Galerkin procedure the DG discretization described above
yields the following matrix equation,
\begin{equation}
G^{(d)}\boldsymbol{\psi}^{(d)}+F^{(d)}\boldsymbol{\psi}^{(d)}+M_{t}\boldsymbol{\psi}^{(d)}-\frac{1}{4\pi}M_{s}\boldsymbol{\varphi}=\mathbf{q}^{(d)},\label{eq:matrixFormLinearTransport}
\end{equation}
where $\boldsymbol{\psi}^{(d)}$and $\boldsymbol{\varphi}$ are the vectors
storing the degree of freedom (DOF) values $\{\psi_{e,m}^{d}\}$ and
$\{\varphi_{e,m}\}$, $e=1\dots N_{\kappa},m=1\dots N_{u}$ . The matrices
$G^{(d)},\,F^{(d)},\,M_{t},\,\text{ and }M_{s}$ are block-diagonal
by mesh element indexing; $\mathbf{q^{(d)}}$ is also a block-diagonal
vector. Indexing DOFs using $m$ and $n$ within a mesh element $\kappa_{e}$,
the individual blocks of $G^{(d)}$, $M_{t}$, and $M_{s}$ are, respectively,
\begin{equation}
\left[G_{e}^{(d)}\right]_{m,n}=-\int_{\kappa_{e}}\left(\boldsymbol{\Omega}_{d}\cdot\nabla_{\mathbf{x}}u_{m}^{e}\right)u_{n}^{e}d\mathbf{x},\label{eq:gMatrix}
\end{equation}
\begin{equation}
[M_{t,e}]_{m,n}=\int_{\kappa_{e}}\sigma_{t}(\mathbf{x})u_{m}^{e}u_{n}^{e}d\mathbf{x}\label{eq:mtMatrix}
\end{equation}

\begin{equation}
[M_{s,e}]_{m,n}=\int_{\kappa_{e}}\sigma_{s}(\mathbf{x})u_{m}^{e}u_{n}^{e}d\mathbf{x}\,\,.\label{eq:msMatrix}
\end{equation}
The source spatial moments vector, $\boldsymbol{q}^{(d)}=\{q_{e,m}^{(d)}\},$
is given by
\begin{equation}
q_{e,m}^{(d)}=\int_{\kappa_{e}}q_{d}(\mathbf{x})u_{m}^{e}d\mathbf{x}  - \int_{\partial \kappa_{e} \cap \partial \mathcal{D}}
  \frac{1}{2} \left( \boldsymbol{\Omega}_{d} \cdot \mathbf{n}_{e} + \left|\boldsymbol{\Omega}_{d}\cdot\mathbf{n}_{e} \right| \right)
                u_{m}^{e} \psi_{d,\text{inc}}^{(d)} d \boldsymbol{s}     \,\,.          \label{eq:qVectorDef}
\end{equation}
The face matrix in (\ref{eq:matrixFormLinearTransport}), $F^{(d)}$,
has both block-diagonal components, $F_{e}^{(d)}$,
\begin{equation}
\left[F_{e}^{(d)}\right]_{m,n}=\int_{\Gamma_{e,e}}\left(\frac{1}{2}\boldsymbol{\Omega}_{d}\cdot\mathbf{n}_{e}+\frac{1}{2}\left|\boldsymbol{\Omega}_{d}\cdot\mathbf{n}_{e}\right|\right)u_{m}^{e}u_{n}^{e} d\boldsymbol{s},\label{eq:diagonalFaceMatrixBlocks}
\end{equation}
as well as off-diagonal block components, $F_{e,e'}^{(d)}$. For two
neighboring mesh elements $\kappa_{e}$ and $\kappa_{e'}$ that share
face $\Gamma_{e,e'}$, $F_{e,e'}^{(d)}$ is defined as
\begin{equation}
\left[F_{e,e'}^{(d)}\right]_{m,n}=\int_{\Gamma_{e,e'}}\left(\frac{1}{2}\left|\boldsymbol{\Omega}_{d}\cdot\mathbf{n}_{e}\right|-\frac{1}{2}\boldsymbol{\Omega}_{d}\cdot\mathbf{n}_{e}\right)u_{m}^{e}u_{n}^{e'}d\boldsymbol{s},\,\,\,\,e'\neq e.\label{eq:offDiagonalFaceMatrixBlocks}
\end{equation}
It is critical to rememember that in (\ref{eq:diagonalFaceMatrixBlocks})
and (\ref{eq:offDiagonalFaceMatrixBlocks}), $\mathbf{n}_{e}$ points
from element $\kappa_{e}$ into $\kappa_{e'}$, \textit{always}. 

The individual blocks defined in (\ref{eq:gMatrix})-(\ref{eq:msMatrix})
can be expressed in terms of integration over the reference element
$\hat{\kappa}$ using the mesh element transformation defined in (\ref{eq:elementTransformation}),
\begin{equation}
\left[G_{e}^{(d)}\right]_{m,n}=-\int_{\hat{\kappa}}\left(\boldsymbol{\Omega}_{d}\cdot J^{-T}(\boldsymbol{\xi})\nabla_{\boldsymbol{\xi}}\hat{u}_{m}\right)\hat{u}_{n}\left|J^{e}(\boldsymbol{\xi})\right|d\boldsymbol{\xi}\,\,,\label{eq:localGMatrix}
\end{equation}

\begin{equation}
\left[M_{t,e}\right]_{m,n}=\int_{\hat{\kappa}}\sigma_{t}\hat{u}_{m}\hat{u}_{n}\left|J^{e}(\boldsymbol{\xi})\right|d\boldsymbol{\xi}\,\,\,,\text{and}\label{eq:localMtMatrix}
\end{equation}

\begin{equation}
\left[M_{s,e}\right]_{m,n}=\int_{\hat{\kappa}}\sigma_{s}\hat{u}_{m}\hat{u}_{n}\left|J^{e}(\boldsymbol{\xi})\right|d\boldsymbol{\xi}\,\,\,\,,\label{eq:localMsMatrix}
\end{equation}
where $J$ is the Jacobian matrix of the coordinate matrix defined
in (\ref{eq:elementTransformation}):
\begin{equation}
J^{e}(\boldsymbol{\xi})=\frac{\partial\mathbf{T}_{e}(\boldsymbol{\xi})}{\partial\boldsymbol{\boldsymbol{\xi}}}=\sum_{m=1}^{N_{v}}\mathbf{x}_{e,m}\left[\nabla\hat{v}_{m}(\boldsymbol{\xi})\right]^{T}.\label{eq:Jacobian}
\end{equation}
Likewise, the $F^{(d)}$ block element can be evaluated on a reference
element,

\begin{equation}
\left[F_{e}^{(d)}\right]_{m,n}=\int_{\hat{\Gamma}}\left(\frac{1}{2}\boldsymbol{\Omega}_{d}\cdot\mathbf{n}_{e}+\frac{1}{2}\left|\boldsymbol{\Omega}_{d}\cdot\mathbf{n}_{e}\right|\right)\hat{u}_{m}\hat{u}_{n'}\left|J^{\Gamma}(\hat{\boldsymbol{s}})\right|d\hat{\boldsymbol{s}},\label{eq:localFDiagonalMatrix}
\end{equation}

\begin{equation}
\left[F_{e,e'}^{(d)}\right]_{m,n}=\int_{\hat{\Gamma}}\left(\frac{1}{2}\left|\boldsymbol{\Omega}_{d}\cdot\mathbf{n}_{e}\right|-\frac{1}{2}\boldsymbol{\Omega}_{d}\cdot\mathbf{n}_{e}\right)\hat{u}_{m}\hat{u}_{n'}\left|J^{\Gamma}(\hat{\boldsymbol{s}})\right|d\hat{\boldsymbol{s}}\,\,.\label{eq:localFOffDiagonalMatrix}
\end{equation}
In (\ref{eq:localFOffDiagonalMatrix}) and \ref{eq:localFDiagonalMatrix},
$n'$ is used in the reference element integration to emphasize that
the reference element indexing in mesh element $\kappa_{e}$ is not
the same as in $\kappa_{e'}$. The face Jacobian $J^{\Gamma}$ is
derived through the $dim-1$ versions of (\ref{eq:elementTransformation})
and (\ref{eq:Jacobian}).

The integrands in (\ref{eq:localGMatrix})-(\ref{eq:localFOffDiagonalMatrix})
are typically evaluated using numerical quadrature. In particular,
Gauss-Legendre quadrature is often used with great success to evaluate
the polynomial (assuming constant material properties) integrands
of (\ref{eq:localGMatrix})-(\ref{eq:localMsMatrix}). However, for
linear mesh elements in 3D, and more generally HO mesh elements in
multi-dimensional geometry, (\ref{eq:localFDiagonalMatrix}) and (\ref{eq:localFOffDiagonalMatrix})
have non-smooth integrands. The $F_{e}^{(d)}$ and $F_{e,e'}^{(d)}$
integrands are non-smooth on re-entrant faces where $\boldsymbol{\Omega}_{d}\cdot\mathbf{n}_{e}$
changes sign. In such cases standard quadrature schemes designed to
integrate smoothly varying functions can converge slowly. To achieve
a simple and rapidly convergent integration scheme, one can use Romberg
integration to rapidly evaluate the piecewise smooth integrands of
(\ref{eq:localFDiagonalMatrix})-(\ref{eq:localFOffDiagonalMatrix})
\cite{numerical_book}.

\section{Graph-Based Transport Sweep Algorithm}

We now describe the manner in which we solve (\ref{eq:matrixFormLinearTransport}).
First, we note that while it is possible to solve the linear system
described by (\ref{eq:matrixFormLinearTransport}), the coupling of
all $\psi^{(d)}$ to all other directions implied by the definition
of (\ref{eq:snScalarFlux}) requires that the angular flux for every
direction be solved simultaneously, a computationally intractable
feat for all but the coarsest of problem resolutions. Rather than
solve for the $N_{\Omega}$ angular fluxes simultaneously, (\ref{eq:matrixFormLinearTransport})
is typically solved using fixed point iteration, 
\begin{equation}
\left[G^{(d)}+F^{(d)}+M_{t}\right]\boldsymbol{\psi}_{(\ell+1)}^{(d)}=\mathbf{q}^{(d)}+\frac{1}{4\pi}M_{s}\boldsymbol{\varphi}_{(\ell)}\label{eq:sourceIteration}
\end{equation}
where for iteration index $\ell$, (\ref{eq:sourceIteration}) is
used to lag the angle coupling,
\[
\boldsymbol{\varphi}_{(\ell)}=\sum_{d=1}^{N_{\Omega}}\omega_{d}\,\boldsymbol{\psi}_{(\ell)}^{(d)}.
\]
Equation (\ref{eq:sourceIteration}) is typically referred to as source
iteration in transport literature \cite{lewis_book}. While it is
more tractable to solve (\ref{eq:sourceIteration}) by forming and
exactly applying the ``streaming + collision operator''
\begin{equation}
\left[G^{(d)}+F^{(d)}+M_{t}\right]\,,\label{eq:streamingCollisionOperator}
\end{equation}
or left hand side of (\ref{eq:sourceIteration}), it is more common
to invert the streaming+collision operator for each angle on a mesh
element-by-mesh element basis in what is reffered to as a transport
sweep. In the following discussion we assume that, for a given direction
$\boldsymbol{\Omega}{}_{d}$, each mesh element index corresponds
to the order in which the unknowns are solved. In 2D, for convex domains
discretized with convex mesh elements, or for Cartesian product meshes
in both 2D and 3D, it is always possible to choose a mesh element
ordering for which the face matrix $F_{e,e'}^{(d)}=0$ whenever $e'>e$.
In this case, the ``streaming + collision operator'' is block lower-triangular.
Its exact inverse can be applied by solving (\ref{eq:matrixFormLinearTransport})
element-by-element, in the order given by the mesh ordering for direction
$\boldsymbol{\Omega}_{d}$.

For general meshes, in particular HO meshes, a mesh element ordering
for which the matrix $F^{(d)}$ is block lower triangular may not
exist. Instead, we find a mesh element ordering for which the block
upper-triangular matrix, $F_{>}^{(d)}$, is as small as possible,
where we have decomposed $F^{(d)}$ in an alternative manner,
\begin{equation}
F^{(d)}=F_{\leq}^{(d)}+F_{>}^{(d)}\,\,.\label{eq:fSplitting}
\end{equation}
Once such a mesh ordering is found, we perform a lagged transport
sweep to update $\psi_{(\ell+1)}^{(d)}$,
\begin{equation}
\left[G^{(d)}+F_{\leq}^{(d)}+M_{t}\right]\boldsymbol{\psi}_{(\ell+1)}^{(d)}=\mathbf{q}^{(d)}+\frac{1}{4\pi}M_{s}\boldsymbol{\varphi}_{(\ell)}-F_{>}^{(d)}\boldsymbol{\psi}_{(\ell)}^{(d)}\,.\label{eq:laggedSourceIteration}
\end{equation}
The iteration scheme described in (\ref{eq:laggedSourceIteration})
is equivalent to solving (\ref{eq:sourceIteration}) using a block
Gauss-Seidel solver on a specific ordering of elements. 

We now pose the problem of choosing a good element ordering in a graph-theoretical
sense. In particular, for each quadrature direction $\boldsymbol{\Omega}_{d}$,
we construct a directed graph whose vertices $V_{e}$ correspond to
mesh elements $\kappa_{e}$. There exists a weighted oriented edge
between mesh element vertices $\left(V_{e},V_{e'}\right)$ in the
graph if $F_{e,e'}^{(d)}$ is nonzero. The choice of edge weight,
$z_{e,e'}$, connecting $(V_{e},V_{e'})$ is critical in determining
the convergence behavior of (\ref{eq:laggedSourceIteration}). 
In this work, we consider three choices for $z_{e,e'}$:
\begin{enumerate}
\item Unity, $z_{e,e'}=1$;
\item \cred{ Face, $z_{e,e'} = \left\Vert F_{e,e'} \right\Vert$; } and
\item SigInvFace,
\begin{equation}
\label{eq:graphWeight}
z_{e,e'}=\left\Vert \left(M_{t,e}\right)^{-1}F_{e,e'}\right\Vert,
\end{equation}
\end{enumerate}
\cred{
where $\| A \| = \max_{i,j} |a_{ij}|$ denotes the maximum element norm of a matrix. }
Equation (\ref{eq:graphWeight}) accounts for the importance of effective
streaming area while also taking into account that radiation is less
likely to stream across optically thick mesh elements.

In order to compute a good element ordering, we solve the corresponding
maximum acyclic subgraph problem, also known as the feedback arc set
(FAS) problem \cite{Eades1993}.
In particular, each element ordering corresponds to
some number of forward edges $\left(V_{e},V_{e'}\right)$, $e'>e$,
and some number of backward edges $\left(V_{e},V_{e'}\right)$, $e'\leq e$.
We find an element ordering that minimizes the weight associated with
the forward edges,
\[
\sum_{e}\sum_{e'<e}z_{e,e'},
\]
because the collection of forward edges corresponds to the upper block-triangular
part of the matrix. This collection defines the outgoing face matrices
$F_{e,e'}$ that are applied to lagged data within the transport sweep
(\ref{eq:laggedSourceIteration}). Although this problem is known
to be NP-hard (see, for example, \cite{Eades1993}), efficient heuristic
methods can be used to compute an approximate solution. 

To determine a mesh-element sweep ordering, we first compute all strongly
connected components (SCCs) of the graph \cite{tarjan}. Then, we
use the fact that a globally minimum FAS can be determined by independently
computing on each SCC a minimal FAS. Since most SCCs involve only
a few elements, this often allows the computation of a minimal FAS
without approximation using integer programming. If the number of
vertices in individual SCCs is larger than a constant, say 10, we
determine an (approximate) minimum FAS using the heuristic algorithm
of \cite{Eades1993}.

\section{Results}
\label{sec_results}
\cred{
In this section we present results in the context of $S_N$ transport on
high-order curved meshes.
We focus on the properties of cycle breaking and iterative convergence with
respect to the mesh curvature and sizes of cycles.
Discussion about an alternative approach to handle curved meshes,
namely, transition to a low-order linear mesh, is also included.

Convergence properties of the method for the cases of scattering-dominated
regimes is beyond the scope of this article.
This subject is addressed in \cite{dsa_sweeps} where we present
a new diffusion synthetic acceleration (DSA) preconditioner
for DG discretizations on high-order meshes. 

Our implementation of the algorithms presented in the previous sections
utilizes the MFEM finite element library \cite{mfem}. }

\subsection{Nested Annuli Test Problem}

\begin{figure}[!ht]
\centering\caption{Annular test problem geometry.\label{fig:annularProblemGeometry}}
\includegraphics[width=6cm]{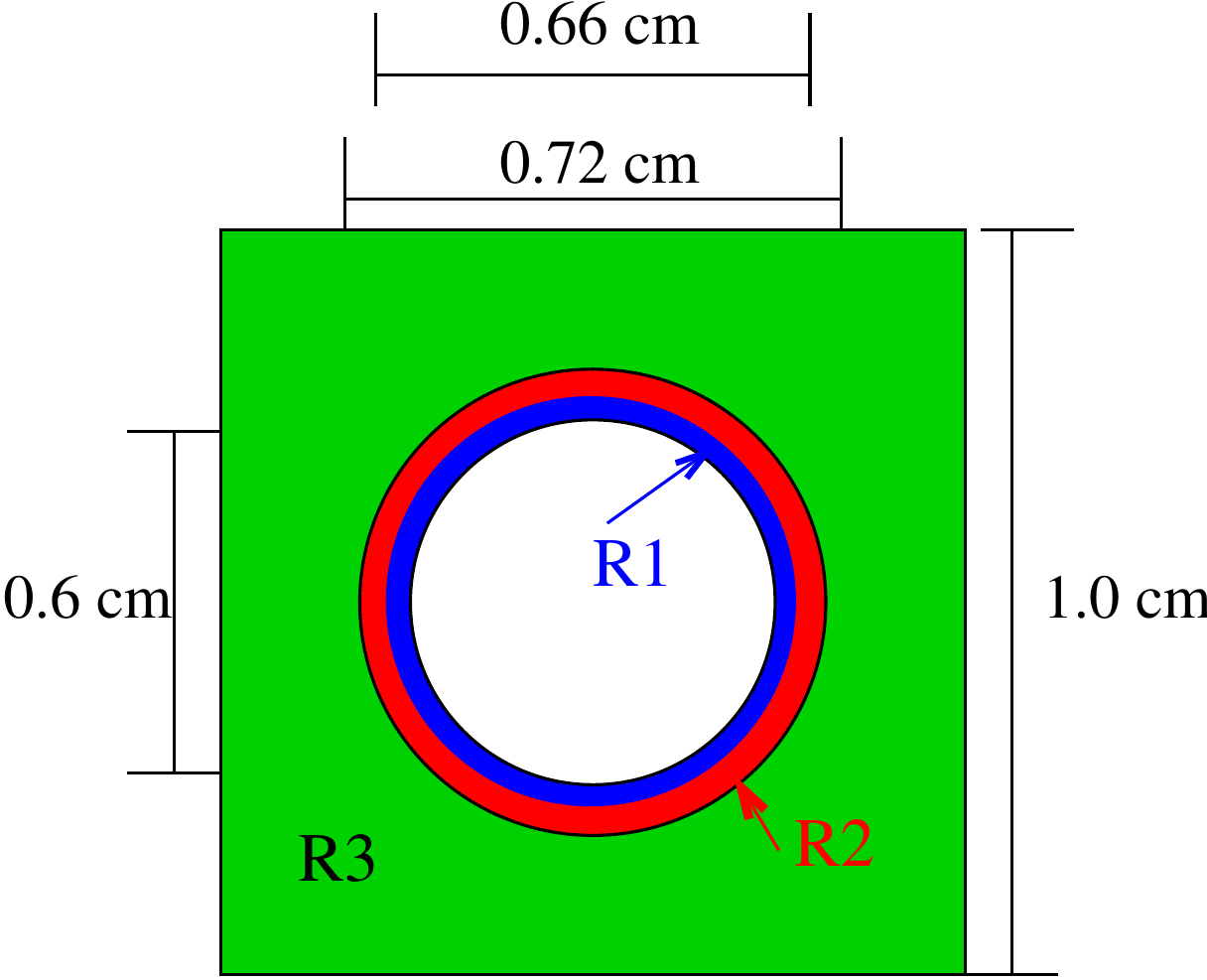}
\end{figure}

\begin{figure}[!ht]
\centering
\caption{Meshes for annular test problem.}
\begin{subfigure}[b]{0.48\textwidth}
    \includegraphics[width=0.95\textwidth,clip=true,trim=3.5in 2.5in 1.75in 1.25in]{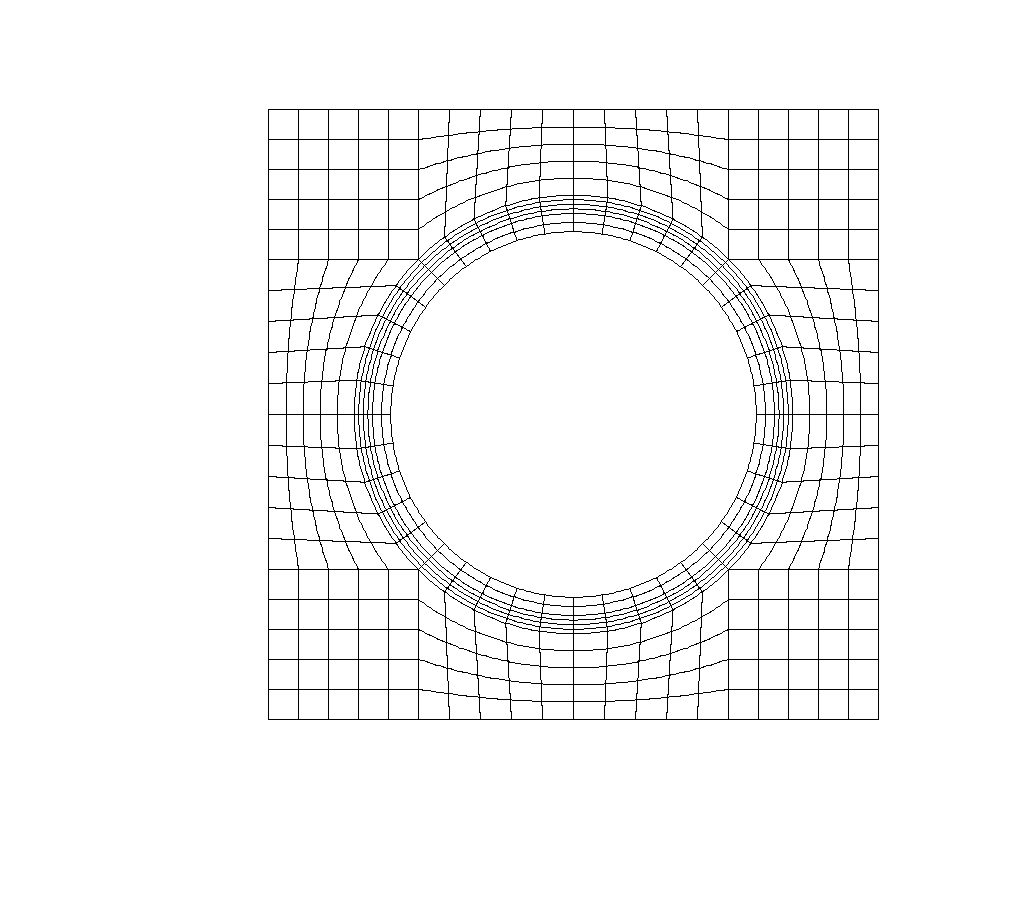}
    \caption{Coarse HO mesh.}
     \label{fig:annulusHoCoarseMesh}
\end{subfigure}
\begin{subfigure}[b]{0.48\textwidth}
    \includegraphics[width=0.95\textwidth,clip=true,trim=3.5in 2.5in 1.75in 1.25in]{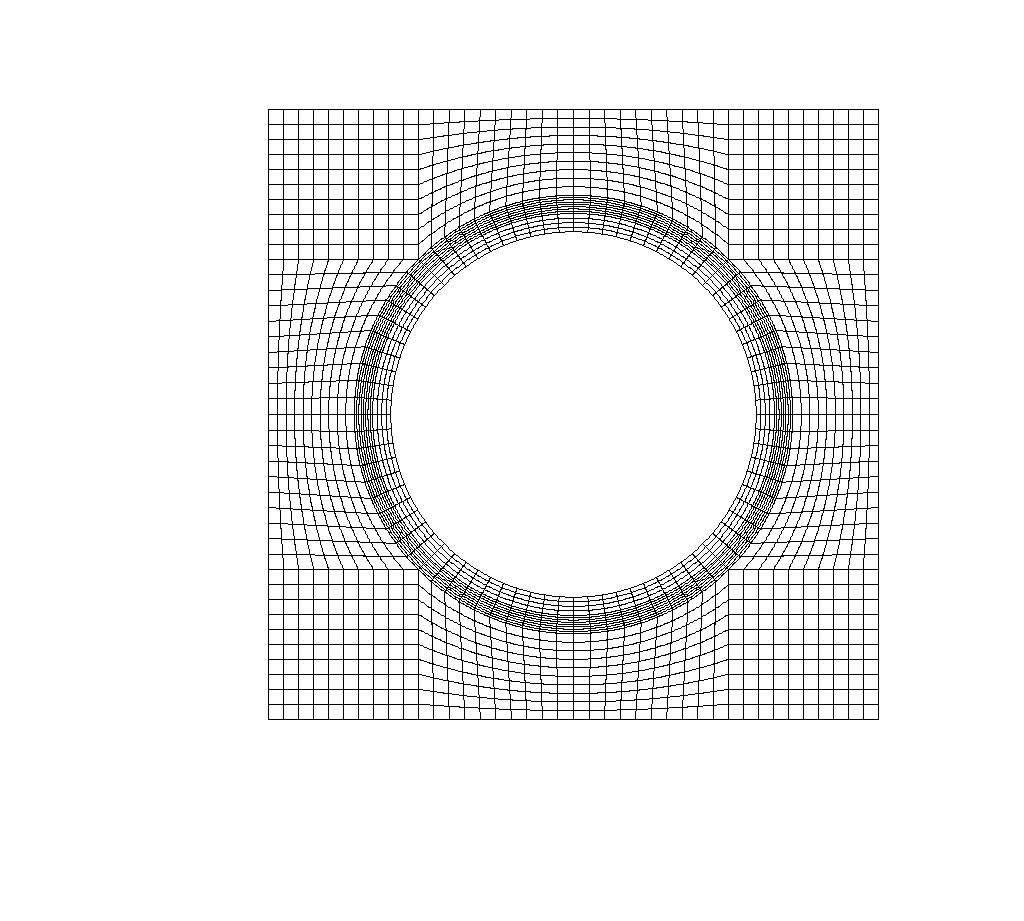}
    \caption{Refined HO mesh.}
     \label{fig:annulusHoRefinedMesh}
\end{subfigure}
\begin{subfigure}[b]{0.48\textwidth}
    \includegraphics[width=0.95\textwidth,clip=true,trim=3.5in 2.5in 1.75in 1.25in]{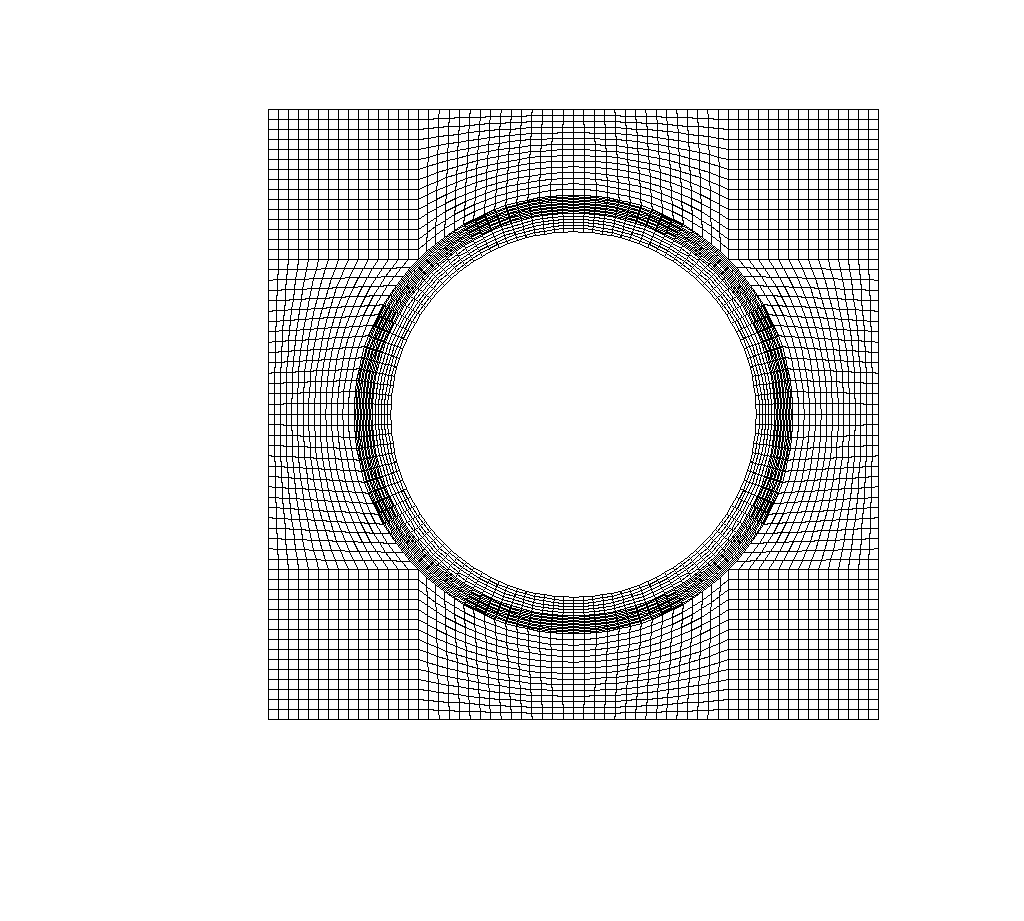}
    \caption{Coarse straightened mesh.}
     \label{fig:annulusStraightenedCoarseMesh}
\end{subfigure}
\begin{subfigure}[b]{0.48\textwidth}
    \includegraphics[width=0.95\textwidth,clip=true,trim=3.5in 2.5in 1.75in 1.25in]{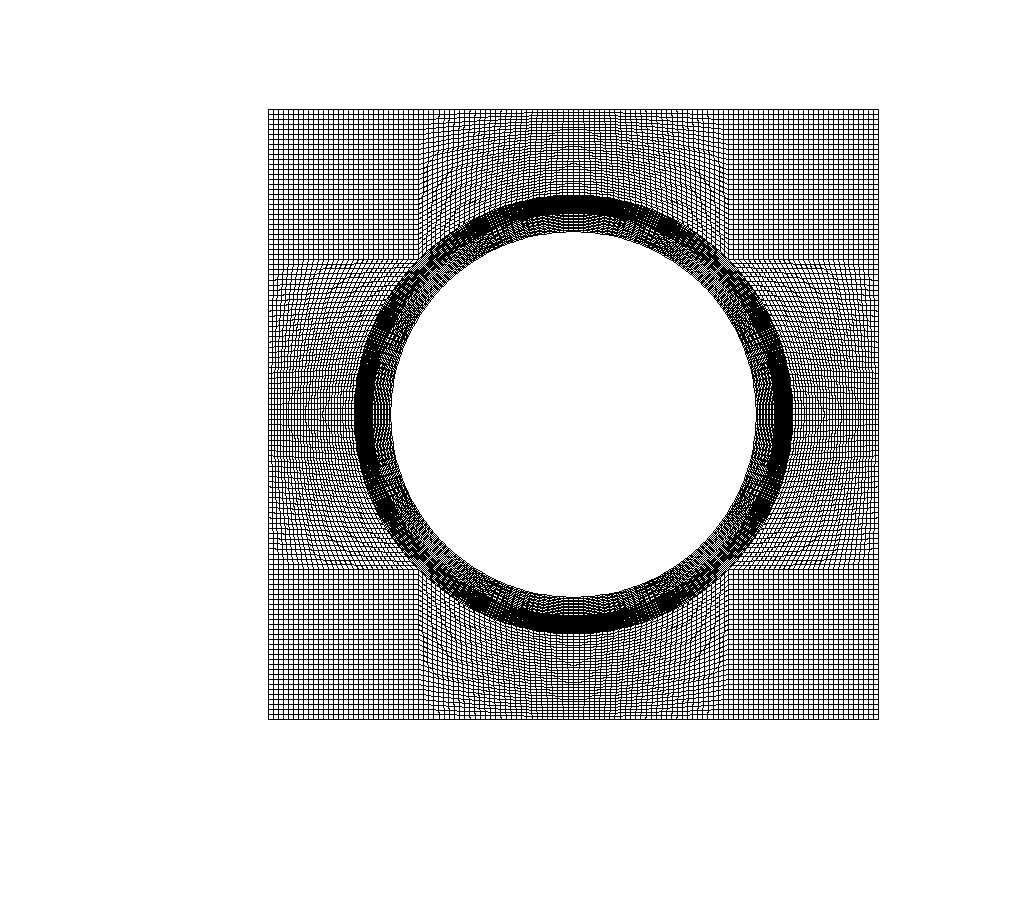}
    \caption{Refined straightened mesh.}
     \label{fig:annulusStraightenedRefinedMesh}
\end{subfigure}
\label{fig:annularMeshes}
\end{figure}

Our first example is a three region test problem consisting of two nested annuli surrounded by a square domain, similar to the geometry of a pressurized water reactor pincell. 
It is important to note that the region interior to the annulus is not included in our problem domain. 
The dimensions of the problem are shown in Fig. \ref{fig:annularProblemGeometry} and material properties are given in Table \ref{tbl:annularMaterialProperties}.
\begin{table}[!ht]
\centering
\caption{Annular test problem material properties.\label{tbl:annularMaterialProperties}}
\begin{tabular}{|c|c|c|c|}
\hline 
Region & R1 & R2 & R3 \\
\hline 
\hline 
$\Sigma_{s}\,\,\left[cm^{-1}\right]$ & 1.0 & 1.0 & 1.0 \\
\hline 
$\Sigma_{t}\,\,\left[cm^{-1}\right]$ & 2.0 & 2.2 & 2.4 \\
\hline 
\end{tabular}
\end{table}
We solve this problem on a series of meshes of cubic elements using cubic DG representations of the angular flux,
\begin{itemize}
\item HO coarse- 540 cubic mesh elements (Fig. \ref{fig:annulusHoCoarseMesh}),
\item HO refined- 2160 cubic mesh elements (Fig. \ref{fig:annulusHoRefinedMesh}),
\end{itemize}
and low-order, spatially refined meshes,
\begin{itemize}
\item Straightened coarse- 4860 linear mesh elements (Fig. \ref{fig:annulusStraightenedCoarseMesh}),
and
\item Straightened refined- 19440 linear elements (Fig. \ref{fig:annulusStraightenedRefinedMesh}).
\end{itemize}

\begin{figure}[ht]
\centering
    \caption { Qualitative description of the straightening procedure.}
    \includegraphics[width=0.85\textwidth]{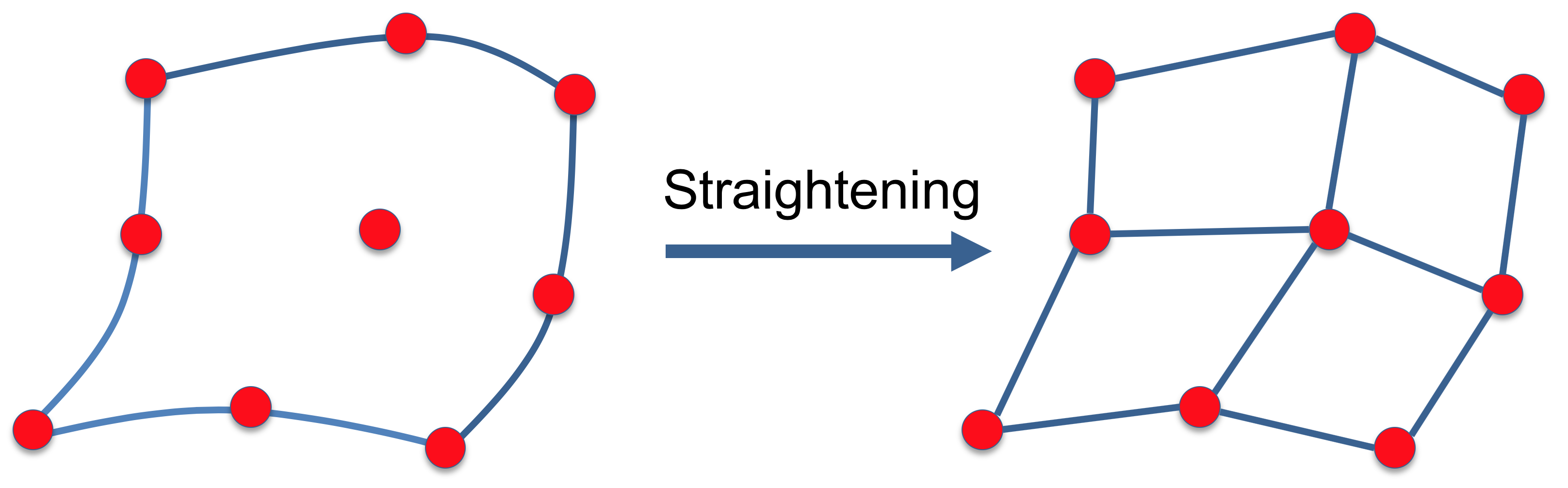}
\label{fig:straightening} 
\end{figure}

The low-order, spatially refined, or ``straightened,'' meshes are generated in
MFEM \cite{mfem} by projecting a HO mesh onto a linear mesh.
This is achieved by dividing each element $N_{ref}$ times in every direction.
The vertices of these sub-divided elements are then taken to
be the nodes of the linear mesh.
This process is demonstrated qualitatively in Fig. \ref{fig:straightening}.

Straightening a HO mesh and solving on the LO mesh using LO finite elements is an alternative to solving
the linear transport equations directly on the HO mesh using HO finite elements. 
Such an approach is appealing in that it,
\begin{enumerate}
\item permits the use of existing, highly-optimized transport solvers designed for meshes of linear elements,
\item likely results in fewer and smaller mesh cycles, if our graph based sweep is too inefficient, and
\item probably generates fewer negative angular flux solutions, because the mesh would be more spatially refined, and
low-order finite elements typically generate fewer negativities in spatially under-resolved regions compared with HO finite elements.
\end{enumerate}
\cred{
The above benefits hold true under several {\em conditions}
that are listed below.
We stress that, up to our knowledge, the current literature doesn't contain
any methods for transition between HO and LO meshes that could provide
all of the listed conditions: }

\begin{enumerate}
\item a spatially refined, \cred{valid} LO mesh can be easily generated from any given HO mesh,
\item mapping materials and interfaces from the HO mesh to the LO mesh does not degrade solution accuracy, 
\item the relative speed of existing, transport solvers operating on a refined LO mesh can compensate for the increase in number of unknowns, and
\item the computational time required for mapping of fields from the HO mesh to the LO mesh and back is trivial.
\end{enumerate}

Using the method of manufactured solutions \cite{lingus,mms}, we first demonstrate that HO solutions on HO meshes are more accurate than LO solutions on LO refined meshes on problems with curved boundaries.
We impose a solution of the form
\cred{
\begin{equation}
\label{eq:psiAnnularMMS}
\psi_{d}\left(x, y\right)=\left[\frac{1}{2}(1+x^{2}+y^{2})+\cos\left(3x+\frac{3y}{2}\right)\right]
  \left[\boldsymbol{\Omega}_{d, \mu}^{2}+
        \boldsymbol{\Omega}_{d, \eta}\right],
\end{equation}
where $\boldsymbol{\Omega}_{d, \mu}$ and $\boldsymbol{\Omega}_{d, \eta}$
are the direction cosines with respect to the $x-$ and $y-$axes, respectively. }
Using $S_{4}$ level-symmetric angular quadrature, we tabulate the $L_{2}$ error
of the scalar flux solution in Table \ref{tbl:annulusErrors}.
As shown in Table \ref{tbl:annulusErrors}, HO transport solutions on HO meshes are orders of magnitude more accurate than LO refined solutions,
despite the HO methods using significantly fewer unknowns.
\begin{table}
\centering
\caption{Comparison of $L_{2}$ errors for MMS on annular test problem.}
\begin{tabular}{|c|c|c|}
\hline 
                                 &      Error               & Number of  \\
{}                               &                             &    $\widetilde{\varphi}(\boldsymbol{x})$ Unknowns  \\
\hline 
Straightened Coarse  & $2.603\times10^{-5}$ &     19940      \\
\hline 
Straightened Refined  & $4.946\times10^{-6}$ &   77760      \\
\hline 
HO Coarse                  & $2.573 1\times10^{-7}$ &   8640      \\
\hline 
HO Refined                 & $1.491\times10^{-8}$ &       34560      \\
\hline 
\end{tabular}
\label{tbl:annulusErrors}
\end{table}

The number of mesh cycles for each HO mesh is summarized in Table \ref{tbl:annuli_mesh_cycles}.
Note that, given the symmetry of this problem and the quadrature set used, we only provide the mesh-cycle information for two angle types, 
\begin{enumerate}
\item the ``diagonal'' ordinates: 
$\boldsymbol{\Omega}_{d,\mu} = \boldsymbol{\Omega}_{d,\eta}$, and
\item ``off-diagonal'' ordinates (every other discrete ordinate direction).
\end{enumerate}
\begin{table}
\centering
\caption{Mesh Cycles in Nested Annuli Test Problem. }
\begin{tabular}{|c|c|c|c|c|}
\hline 
     Mesh     & Angle &      Total Cycles            &  Elements  &    Total Edges \\
     {}            &           &   (Simple+{\bf Large})         &    in Large SCC            &      Removed \\
\hline 
\hline 
HO Coarse &   Diagonal        &  12 + {\bf 0}        &     -       &     12 \\
      {}           &   Off-Diagonal  &   2 + {\bf 2}        &    [7,7]    &     14 \\
\hline 
\hline 
HO Fine    &   Diagonal        &  24 + {\bf 0}        &     -       &     24 \\
      {}         &   Off-Diagonal  &   6 + {\bf 4}        &    [3,3,13,13]    &     34 \\
\hline 
\end{tabular}
\label{tbl:annuli_mesh_cycles}
\end{table}
Note that for this particular problem, unity weighting of inter-element dependencies was used, but the choice of graph-edge weighting did not affect the number of edges removed to make an acyclic sweep ordering.

\begin{figure}[!ht]
\centering
\caption{Iterative convergence for the nested annuli test problem.}
\includegraphics[width=0.7\textwidth]{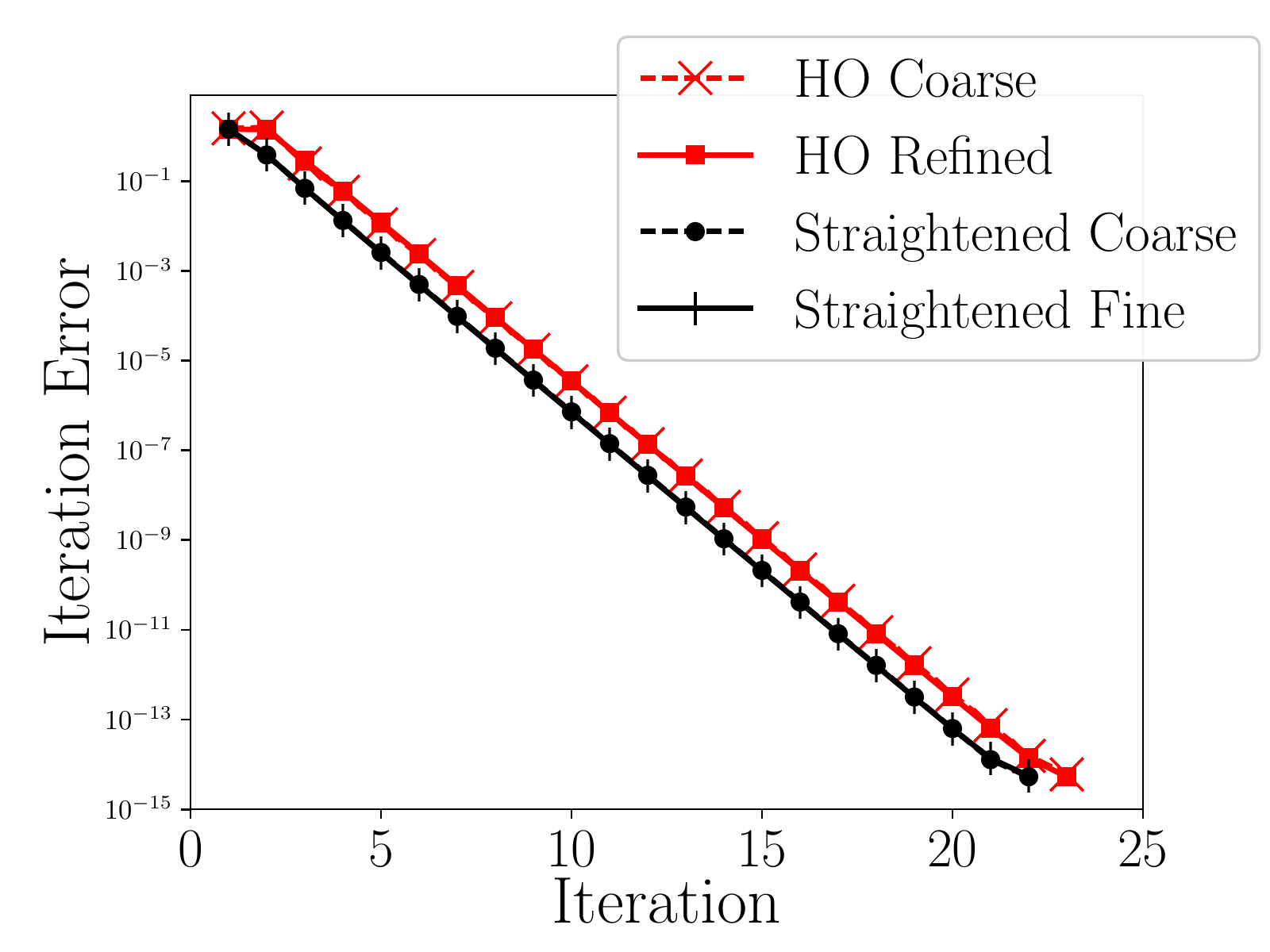}
\label{fig:convergence2Dannulus}
\end{figure}

In Fig. \ref{fig:convergence2Dannulus}, we demonstrate that, despite the HO annular meshes having mesh cycles, with the graph-based sweep ordering, the mesh cycles do not significantly degrade iterative convergence, requiring one additional iteration to reach the same convergence criterion.
Figure \ref{fig:convergence2Dannulus} plots the computed iteration error, taken to be the $L_{\infty}$ norm
of the difference of successive iterates,
\[
\log_{10}\left(\max_{d}\left\Vert \boldsymbol{\psi}_{j+1}^{(d)}-\boldsymbol{\psi}_{j+1}^{(d)}\right\Vert _{\infty}\right)\,\,,
\]
versus the fixed-point iteration index. 
Source iteration for this problem is terminated when the norm of the difference of successive iterates is near machine precision, $10^{-14}$. 

The nested annuli test problems shows that HO transport on HO meshes can be more accurate with fewer unknowns that LO transport discretization on spatially refined grids, as shown by Table \ref{tbl:annulusErrors}.
Additionally, the nested annuli test problem shows that HO meshes introduce mesh cycles that LO meshes do not possess (see Table \ref{tbl:annuli_mesh_cycles}). However, as shown in Fig. \ref{fig:convergence2Dannulus}, this mesh cycling does not significantly degrade iterative convergence relative to cycle-free, straightened meshes of the same problem.

\subsection{Triple-point Problem}
\label{subsec:2D-triple-point-problem}

\begin{figure}[!ht]
\centering
\caption{HO triple-point meshes.}
\begin{subfigure}[b]{0.49\textwidth}
\centering
\includegraphics[width=0.99\textwidth]{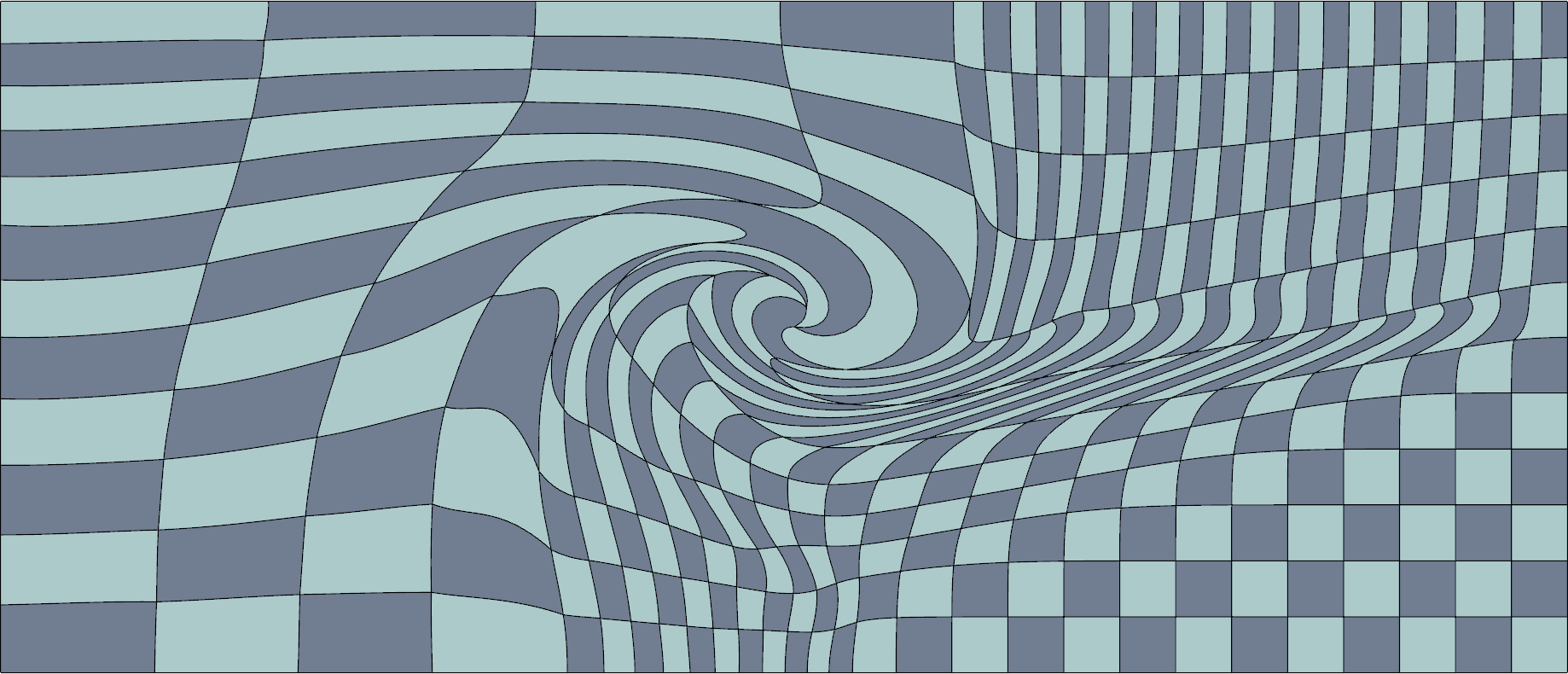}
\caption{Cubic mesh with 336 elements.}
\label{fig:lowResTriplePointMesh}
\end{subfigure}
\begin{subfigure}[b]{0.49\textwidth}
\centering
\includegraphics[width=0.99\textwidth]{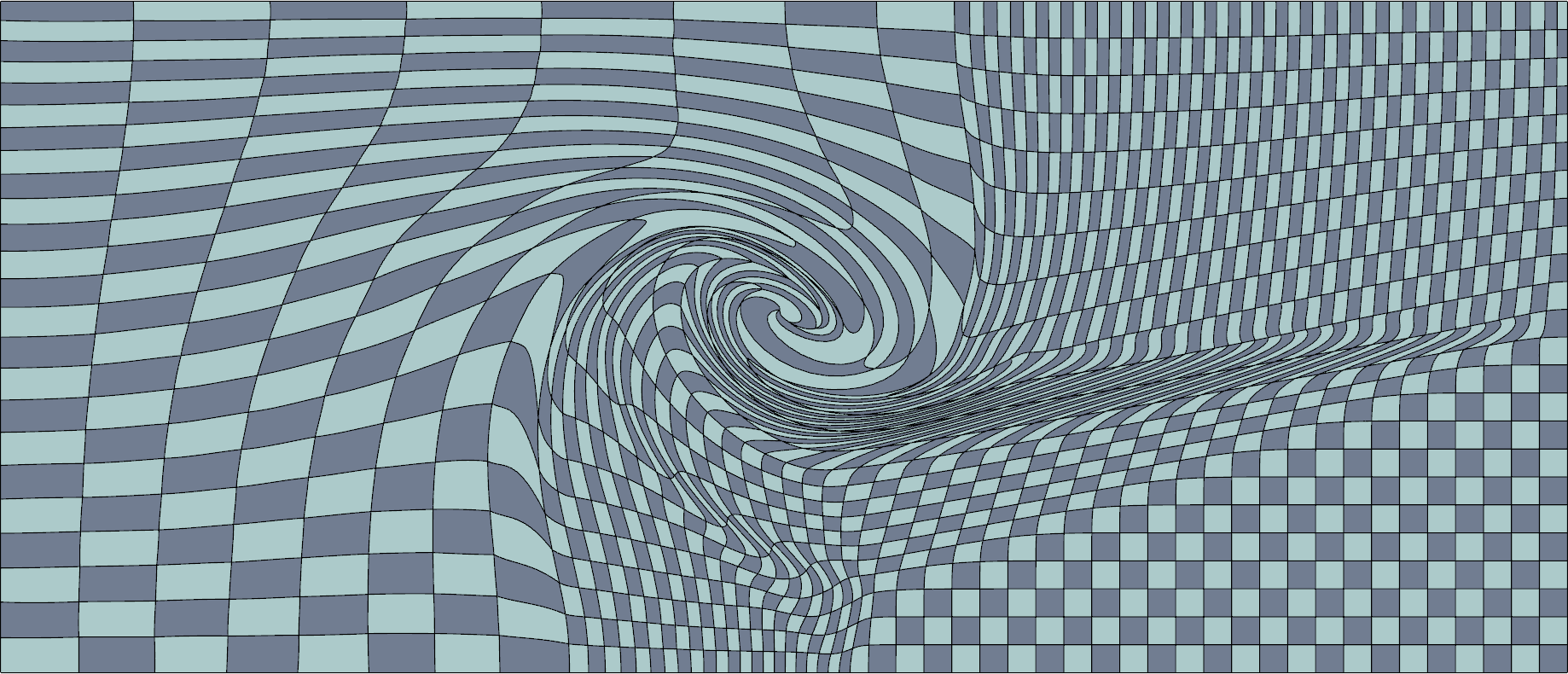}
\caption{Cubic mesh with 1344 elements.}
\label{fig:highResTriplePointMesh}
\end{subfigure}
\label{fig:triple-point-meshes} 
\end{figure}

Our next example uses an underling mesh generated from a purely Lagrangian simulation of the ``triple point'' problem \cite{triplePoint} using cubic discretizations of the kinematic variables, creating third-order meshes.
 The low resolution and high resolution variants of this mesh are shown in Fig. \ref{fig:lowResTriplePointMesh} and Fig.  \ref{fig:highResTriplePointMesh}, respectively.
The scattering and total cross sections are constant, $\Sigma_{s}=1~[cm^{-1}]$ and $\Sigma_{t}=2~[cm^{-1}]$, the problem is driven by an isotropic source $q_{d}\left(\mathbf{x}\right)=1+\sin^{2}\left(2x+y\right)$ and a constant inflow boundary condition of unity for all discrete ordinate directions. 
An $S_{2}$ level symmetric angular quadrature set is used to discretize the problem in angle.

As seen in Fig. \ref{fig:triple-point-meshes}, this problem effectively represents a worst case scenario for mesh cycling in a two-dimensional problem.
We use this problem to demonstrate:
\begin{enumerate}
\item it is possible to sweep on HO meshes with very large strongly connected components,
\item different graph weightings can result in different source iteration convergence rates, and
\item mesh straightening is in general impractical.
\end{enumerate}
First, we give the number of mesh cycles present in the HO meshes we consider in Table \ref{tab:numberSCCs}. 
In Table \ref{tab:numberSCCs} and Table \ref{tbl:edgesRemoved}, $a=\frac{\sqrt{3}}{3}$.
While the number of elements in each SCC is independent of the edge weighting
used, the number of edges removed from the graph is not.
Table \ref{tbl:edgesRemoved} shows the number of edges removed in solving the
FAS sweep ordering problem for each edge-weighting method considered.

\begin{table}[!ht]
\centering
\caption{SCC information for the HO triple-point problem.}
\begin{tabular}{|c|c|c|c|}
\hline 
     Mesh     & Angle &      Total Cycles            &     Elements        \\
     {}            &           &   (Simple+{\bf Large})  &    in Large SCC  \\
\hline 
336 Elements &   $\langle a, a \rangle$        &   0 + {\bf 1}        &     [97]       \\
      {}              &    $\langle -a, a \rangle$        &   1 + {\bf 1}        &    [46]  \\
      {}              &    $\langle a, -a \rangle$        &   1 + {\bf 1}        &    [46]  \\
      {}              &    $\langle -a, -a \rangle$        &   0 + {\bf 1}        &    [97]   \\
\hline 
1344 Elements      &   $\langle a, a \rangle$        &   0 + {\bf 1}        &     [475]       \\
      {}              &    $\langle -a, a \rangle$        &   1 + {\bf 2}        &    [3,184]  \\
      {}              &    $\langle -a, a \rangle$        &   1 + {\bf 2}        &    [3,184]  \\
      {}              &    $\langle -a, -a \rangle$        &   0 + {\bf 1}        &    [475]   \\
 \hline 
\end{tabular}
\label{tab:numberSCCs}
\end{table}

\begin{table}[!ht]
\centering
\caption{Edges removed to define a sweep ordering for the HO triple-point problem.}
\begin{tabular}{|c|c|c|c|c|}
\hline 
     Mesh     & Angle &          \multicolumn{3}{|c|}{ Edge Weighting } \\
     {}            &           &   Unity & Face &  SigmaInvFace   \\
\hline 
336 Elements &   $\langle a, a \rangle$        &  60    &  63  & 72   \\
      {}              &    $\langle -a, a \rangle$       &  25    &  30  & 39   \\
      {}              &    $\langle a, -a \rangle$       &  26    &  31  &  36  \\
      {}              &    $\langle -a, -a \rangle$      &  66    &  64  &  82  \\
      \hline 
1344 Elements&   $\langle a, a \rangle$        &  144  &  170  & 292   \\
      {}              &    $\langle -a, a \rangle$       &  68    &   91    & 143   \\
      {}              &    $\langle a, -a \rangle$       &  81    &   100  &  112  \\
      {}              &    $\langle -a, -a \rangle$      &  202  &   209  &  285  \\
       \hline 
\end{tabular}
\label{tbl:edgesRemoved}
\end{table}

\begin{figure}[!ht]
\centering
\caption{Iterative convergence for triple-point problem.}
\begin{subfigure}[b]{0.49\textwidth}
\centering
\includegraphics[width=0.96\textwidth]{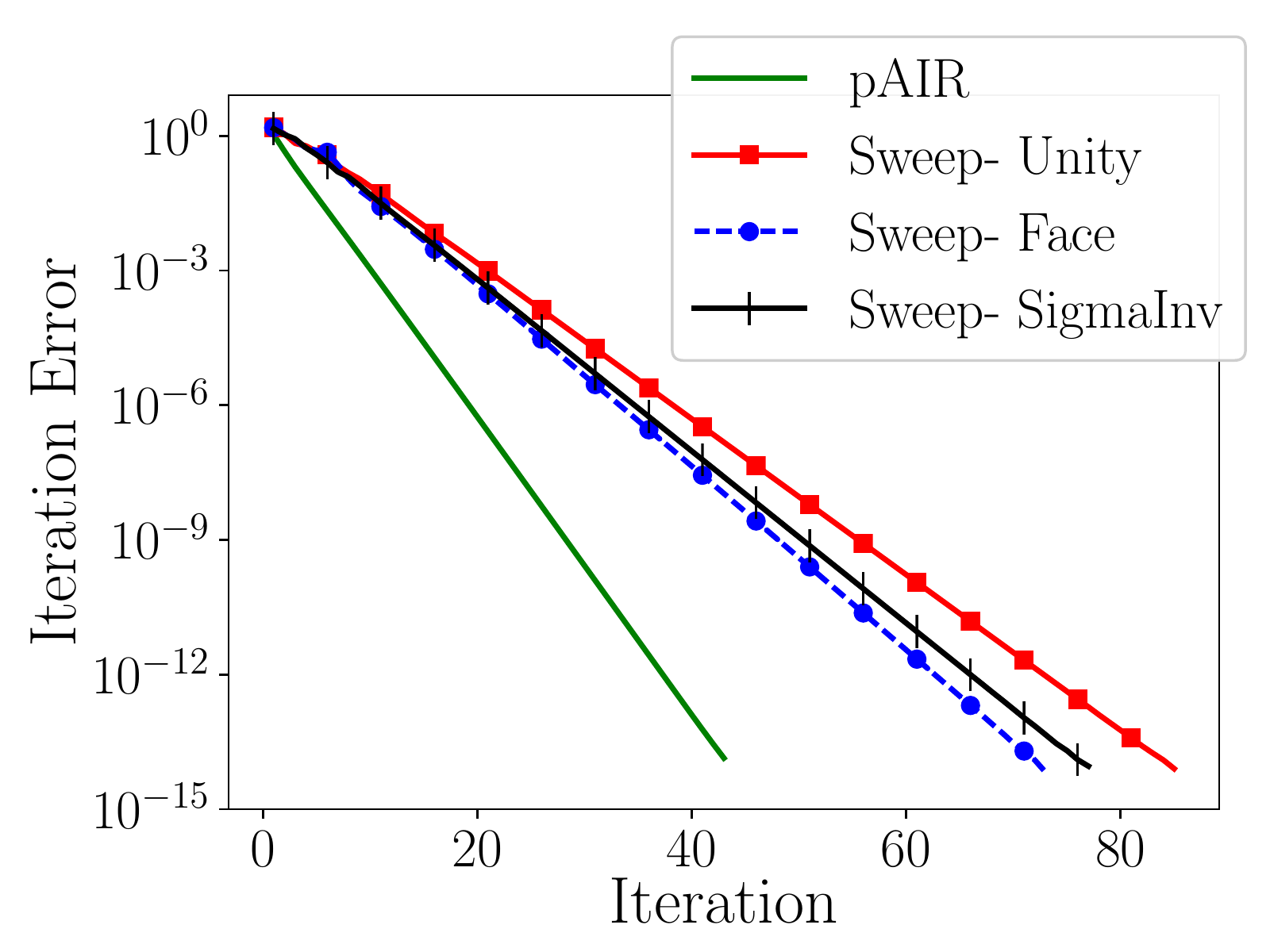}
\caption{336 Element Mesh.}
\end{subfigure}
\begin{subfigure}[b]{0.49\textwidth}
\includegraphics[width=0.96\textwidth]{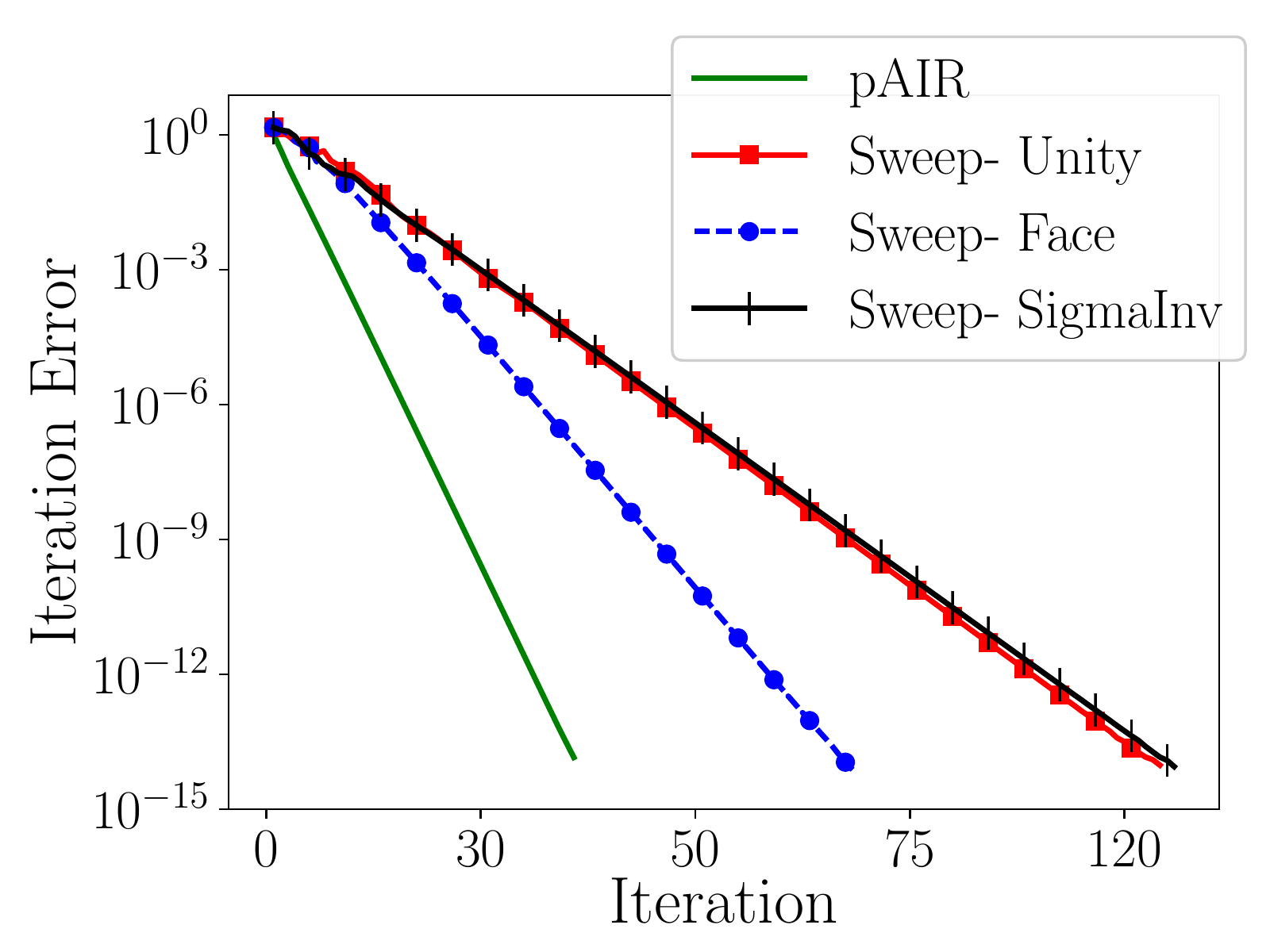}
\caption{1344 Element Mesh.}
\end{subfigure}
\label{fig:triplePointIterativeConvergence}
\end{figure}

To assess how the upwind dependency lagging required to sweep mesh element by element affects overall scattering source iteration, we compare against an ``ideal'' fixed-point, solving Eq. (\ref{eq:sourceIteration}) directly without the Gauss-Seidel iteration introduced by sweeping on HO grids [see Eq. \eqref{eq:laggedSourceIteration}].
We solve Eq. (\ref{eq:sourceIteration}) using a recently developed algebraic multigrid method for non-symmetric matrices, referred to as $p$AIR \cite{air2,air1}.
Iterative convergence for both meshes is plotted in Fig. \ref{fig:triplePointIterativeConvergence}. 
As shown in Fig. \ref{fig:triplePointIterativeConvergence}, sweeping element-by-element does incur an iteration effectiveness penalty, but still converges at a reasonable rate, even for problems with large SCCs.

Finally, we use the HO triple-point meshes to demonstrate that in general,
it is at best impractical to straighten HO meshes.
For this study, we simply seek the minimum refinement level necessary to yield
a low-order, spatially refined mesh that satisfies
$\left \lvert J(\boldsymbol{\xi}) \right \rvert > 0$ at every node.

Straightening of the purely Lagrangian triple-point meshes to yield
positive $\left \lvert J(\boldsymbol{\xi}) \right \rvert$ everywhere in the low-order mesh
is not always possible, e.g., as in the case with the 1344 element mesh.
Furthermore, straightening is likely to be impractical as shown by the 336
element mesh needing 24 levels of refinement to yield strictly positive
$\left \lvert  J(\boldsymbol{\xi}) \right \rvert$.
In particular, refining the HO mesh by a factor of 24, to create a low-order
mesh for low-order transport, would increase the number of transport unknowns
by a factor of 144 for this problem.
Attempts were made to use up to 180 refinements to straighten the 1344 element
HO mesh, however we could not verify that even at that level of refinement the
straightened mesh had \cred{ $|J(\boldsymbol{\xi})| > 0$ at every node. }

\subsection{3D Rayleigh-Taylor problem}
\label{sec_3DRT}
\begin{figure}[!ht]
  \centerline
  {
    \includegraphics[width=0.85\textwidth]{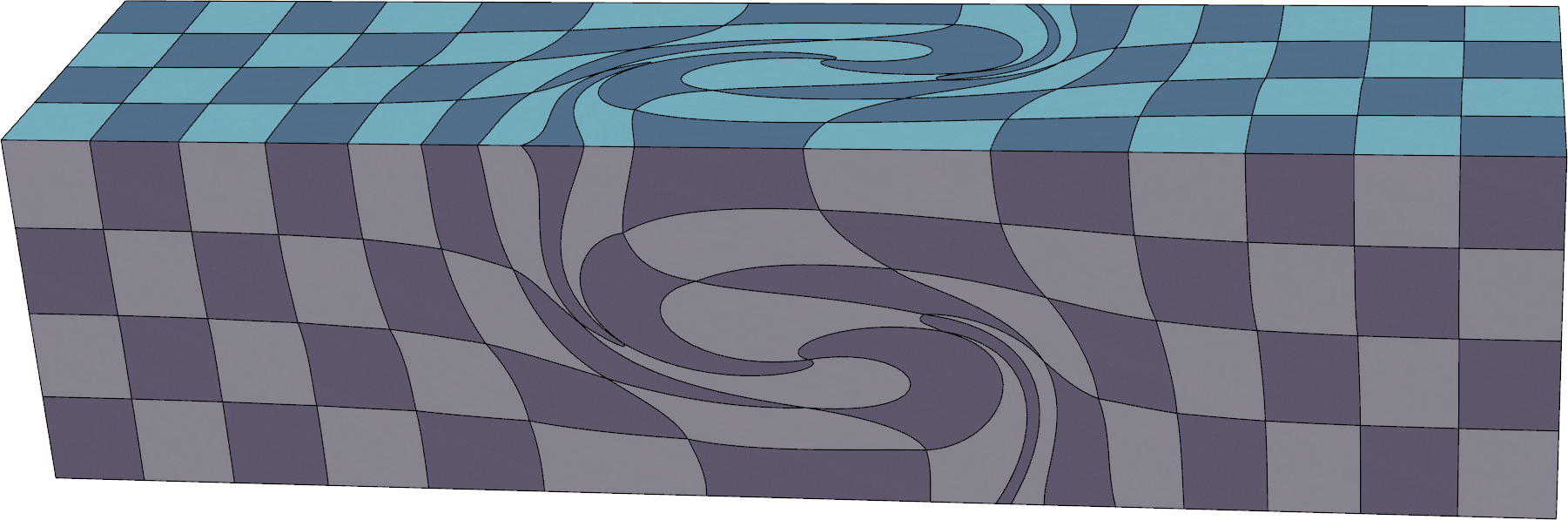}
  } \vspace{4mm}
  \centerline
  {
    \includegraphics[width=0.85\textwidth]{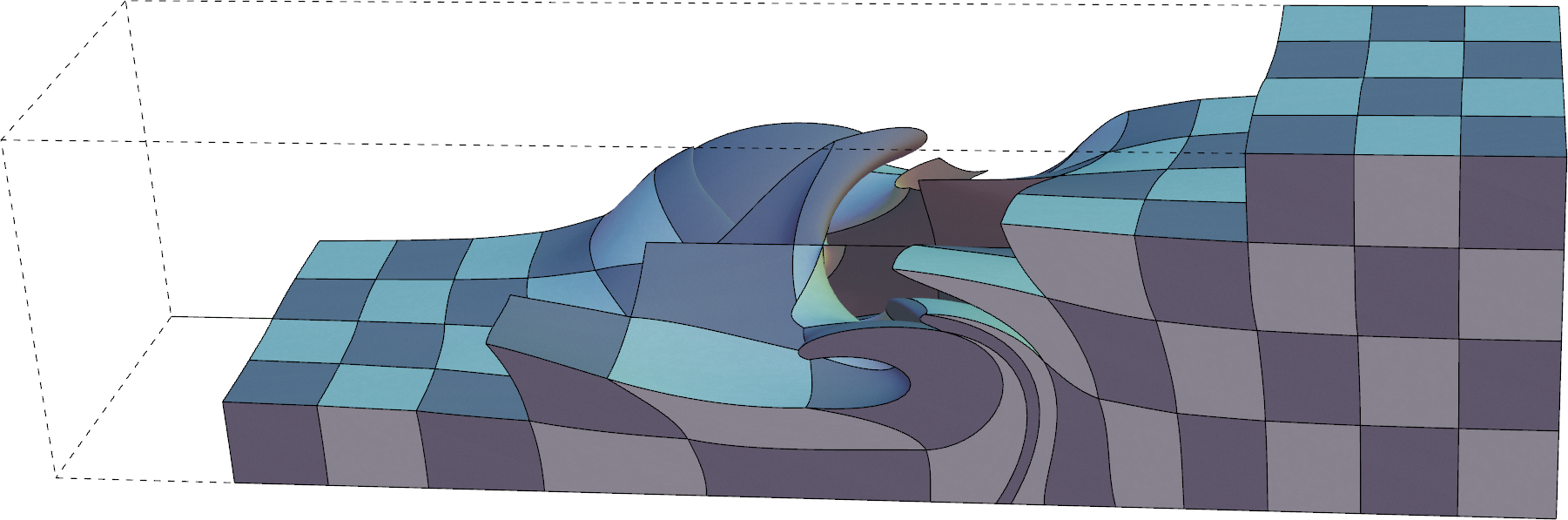}
  }
  \caption
  {
    Cubic 3D mesh (top) resulting from a Lagrangian simulation of the
    Rayleigh-Taylor problem. A subset of the mesh elements is shown on bottom.
  }
  \label{fig_RTmesh} 
\end{figure}

In our final example we consider a three-dimensional third-order mesh generated
via a purely Lagrangian simulation of a 3D variation of the classic
Rayleigh-Taylor instability problem \cite{Dobrev2012}.
The mesh contains 256 zones and is shown in Figure \ref{fig_RTmesh}.
This problem demonstrates that sweeping is still possible on highly-distorted
3D meshes.

Like the triple-point problem, the 3D Rayleigh-Taylor problem also illustrates a
case when straightening the HO mesh is not a feasible alternative, due to the
high curvature of the original mesh.
In particular, obtaining a valid, LO refined version of a coarser variant of this problem, 32 third order zones,
requires refining each HO mesh element by a factor of 15, resulting in a linear mesh containing 108000 elements with 864000 transport spatial degrees of freedom versus 32 HO elements and 2048 transport spatial degrees of freedom.

The scattering and total opacities are constant throughout the domain,
$\sigma_{s}=1$ and $\sigma_{t}=2$, and we set a constant inflow boundary
condition of unity for all discrete ordinate directions.
\cred{
An $S_{2}$ level symmetric angular quadrature set is used to discretize the
problem in angle. }
The problem is driven by the following volumetric source:
\[
q_{d}\left(\mathbf{x}\right) = 1 + \sin^{2}\left( 4x + 2yz + 2y \right) \,.
\]

\begin{figure}[!ht]
  \centering
  \begin{subfigure}[b]{0.45\textwidth}
  \includegraphics[width=0.95\textwidth]{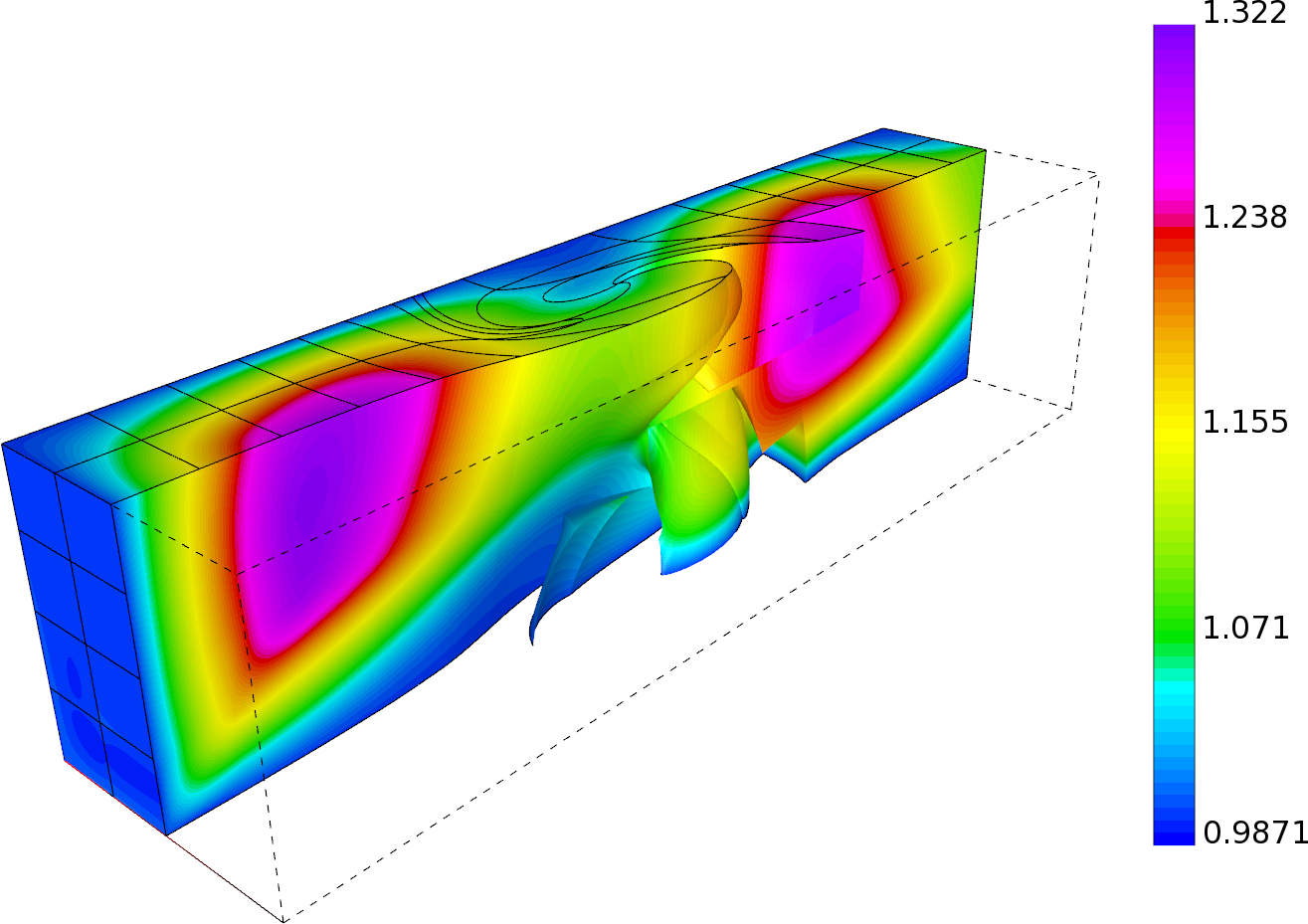} 
  \caption{HO Transport solution.}
  \end{subfigure}
\begin{subfigure}[b]{0.45\textwidth}
    \includegraphics[width=0.95\textwidth]{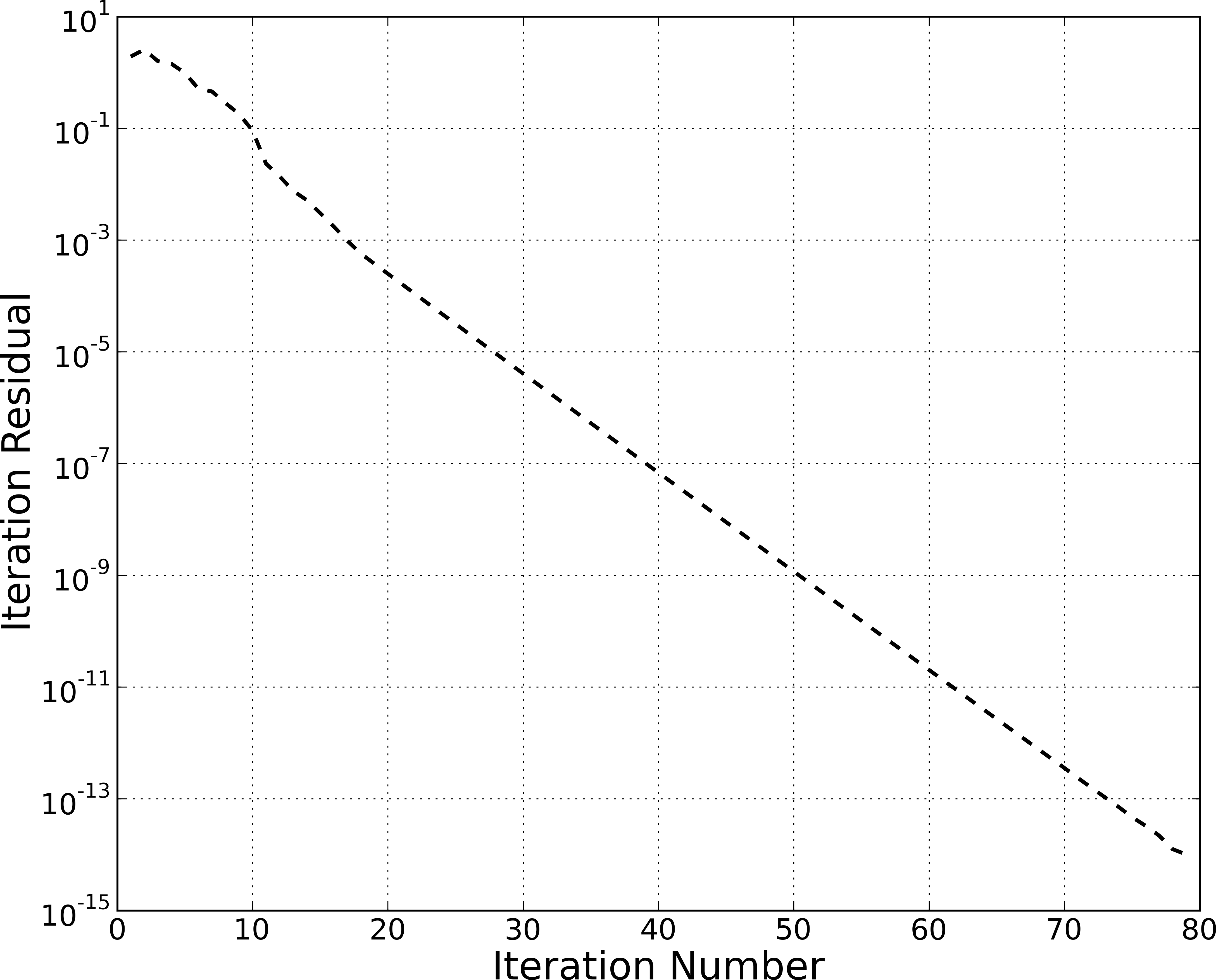}
    \caption{Convergence of iteration error.}
 \end{subfigure}
  \caption{Final solution and convergence of iterative error for sweeping on the Rayleigh-Taylor mesh. }
  \label{fig_RTsltn}
\end{figure}

Each of the eight discrete ordinate direction develops one SCC
containing 79 or 83 elements, out of 256 elements in the full mesh.
\cred{
The number of removed graph edges for each of the three weighting strategies are
listed in Table \ref{tab:edgesRT}.
The Unity, Face, and SigmaInvFace strategies converge in 79, 78, and 83
iterations, respectively, with very similar convergence histories.
The final solution and convergence history with the SigmInvFace strategy
are shown in Figure \ref{fig_RTsltn}. }

\begin{table}[!ht]
\centering
\caption{\cred{Edges removed to define a sweep ordering for
         the 3D Rayleigh-Taylor problem.}}
\begin{tabular}{|c|c|c|c|c|}
\hline 
 Angle & Elements     & \multicolumn{3}{|c|}{ Edge Weighting } \\
       & in Large SCC & Unity & Face &  SigmaInvFace   \\
\hline 
 $\langle a, a, a \rangle$     &  83  &  106   &  104  &  123  \\
 $\langle -a, a, a \rangle$    &  79  &   90   &   93  &  103  \\
 $\langle a, -a, a \rangle$    &  83  &  102   &  103  &  120  \\
 $\langle -a, -a, a \rangle$   &  79  &   94   &   89  &  112  \\
 $\langle a, a, -a \rangle$    &  79  &   95   &   93  &  108  \\
 $\langle -a, a, -a \rangle$   &  83  &  111   &  104  &  123  \\
 $\langle a, -a, -a \rangle$   &  79  &   98   &   89  &  113  \\
 $\langle -a, -a, -a \rangle$  &  83  &  102   &  103  &  120  \\
\hline
\end{tabular}
\label{tab:edgesRT}
\end{table}

\section{Conclusion and Future Work}
\label{sec_future}
Solving the linear transport equations directly on HO grids is a requirement for next-generation radiation hydrodynamics
simulations. \cred{We have demonstrated that, despite a potentially large number of mesh cycles, the linear transport equations can
be solved directly on HO meshes by introducing an appropriate mesh ordering on which to do a transport sweep. The
ordering is determined by approximately solving the feedback arc set (FAS) problem, for some measure of distance between
mesh elements. Three such distances, or ``weightings,'' were explored, the unity weighting, face weighting, and sigma-inverse face weighting
(see Section \ref{sec_3DRT}). In all tested cases, the face weighting was the most effective, in one case offering a near $2\times$
improvement in total iteration count over the unity weighting and sigma-inverse face weighting (see Figure \ref{fig:triplePointIterativeConvergence}).
Solving directly on HO meshes eliminates the need to transition to a straightened linear mesh, improves solution accuracy, and requires
fewer degrees of freedom.}

A future area of research will be the preservation of angular flux positivity when using high-order discontinuous Galerkin spatial discretizations of the linear transport equations in spatially under resolved regions, \cred{as well as a more rigorous understanding of the optimal weighting for the FAS problem.}
We will also explore the challenges HO meshes may pose to diffusion-based synthetic acceleration techniques required for scattering-dominated regimes.
Such capabilities are required by multi-material problems involving optically thick and thin transitions.
The long-term goal is to derive a time-dependent, nonlinear TRT solver utilizing the iterative techniques described within this work.

\section*{Acknowledgements}
The work of T. S. Haut, P. G. Maginot, V. Z. Tomov, T. A. Brunner, and T. S. Bailey was performed under the auspices of the U.S. Department of Energy by Lawrence Livermore National Laboratory under contract DE-AC52-07NA27344 and supported by the LLNL-LDRD Program under Project No. 18-ERD-002.

The work of B. S. Southworth was performed under the auspices of the U.S. Department of Energy by Lawrence Livermore National Laboratory under contract DE-AC52-07NA27344 and supported by sub-contracts B614452, and B627942 of Lawrence Livermore National Security, LLC. 
Additional funding for B. S. Southworth was provided under grant number (NNSA) DE-NA0002376.

This work has been reviewed for unlimited public release as LLNL-JRNL-759880-DRAFT.
  
\section*{Disclaimer}
 This document was prepared as an account of work sponsored by an agency of the United States government. Neither the United States government nor Lawrence Livermore National Security, LLC, nor any of their employees makes any warranty, expressed or implied, or assumes any legal liability or responsibility for the accuracy, completeness, or usefulness of any information, apparatus, product, or process disclosed, or represents that its use would not infringe privately owned rights. Reference herein to any specific commercial product, process, or service by trade name, trademark, manufacturer, or otherwise does not necessarily constitute or imply its endorsement, recommendation, or favoring by the United States government or Lawrence Livermore National Security, LLC. The views and opinions of authors expressed herein do not necessarily state or reflect those of the United States government or Lawrence Livermore National Security, LLC, and shall not be used for advertising or product endorsement purposes.

\bibliographystyle{ans_js}                                                                          \bibliography{pgm_references}

\begin{thebibliography}{10}

\bibitem{mfem}
{MFEM}: Modular finite element methods library.
\newblock {\em http://mfem.org}.

\bibitem{adams_nowak}
M.~L. Adams and P.~F. Nowak.
\newblock Asymptotic analysis of a computational method for time- and
  frequency- dependent radiative transfer.
\newblock {\em Journal of Computational Physics}, 46:366--403, 1998.

\bibitem{Dobrev2018}
Robert~W. Anderson, Veselin~A. Dobrev, Tzanio~V. Kolev, Robert~N. Rieben, and
  Vladimir~Z. Tomov.
\newblock High-order multi-material {ALE} hydrodynamics.
\newblock {\em SIAM J. Sci. Comp.}, 40(1):B32--B58, 2018.

\bibitem{bell_glasstone}
G.~Bell and S.~Glasstone.
\newblock {\em Nuclear Reactor Theory}.
\newblock Van Nostrand Rienhold, Inc., New York, NY, 1970.

\bibitem{Boscheri2016}
Walter Boscheri and Michael Dumbser.
\newblock High order accurate direct arbitrary-{L}agrangian-{E}ulerian
  {ADER}-{WENO} finite volume schemes on moving curvilinear unstructured
  meshes.
\newblock {\em Comput. Fluids}, 136:48--66, 2016.

\bibitem{supinski2017}
B.~R. de~Supinski.
\newblock Sierral the {LLNL} {IBM} {CORAL} system.
\newblock In {\em 12th International Conference on Parallel Processing and
  Applied Mathematics}, September 2017.

\bibitem{Dobrev2013}
V.~Dobrev, T.~Ellis, Tz. Kolev, and R.~Rieben.
\newblock High-order curvilinear finite elements for axisymmetric {L}agrangian
  hydrodynamics.
\newblock {\em Comput. Fluids}, pages 58--69, 2013.

\bibitem{Dobrev2012}
V.~Dobrev, Tz. Kolev, and R.~Rieben.
\newblock High-order curvilinear finite element methods for {L}agrangian
  hydrodynamics.
\newblock {\em SIAM J. Sci. Comp.}, 34(5):606--641, 2012.

\bibitem{Eades1993}
Peter Eades, Xuemin Lin, and W.~F. Smyth.
\newblock {A fast and effective heuristic for the feedback arc set problem}.
\newblock {\em Information Processing Letters}, 47(6):319--323, 1993.

\bibitem{kull}
J.~A.~Rathkopf et~al.
\newblock {KULL}: {LLNL}'s {ASCI} inertial confinement fusion simulation code.
\newblock In {\em {PHYSOR} 2000 American Nuclear Society Topical Meeting on
  Advances in Reactor Physics and Mathematics and Computation into the Next
  Millenium}, May 7-11 2000.

\bibitem{Tomov2016}
Jean-Luc Guermond, Bojan Popov, and Vladimir Tomov.
\newblock Entropy-viscosity method for the single material {E}uler equations in
  {L}agrangian frame.
\newblock {\em Comput. Methods Appl. Mech. Eng.}, 300:402--426, 2016.

\bibitem{triplePoint}
M.~Kucharik, R.~V. Garimella, S.~P. Schofield, and M.~J. Shashkov.
\newblock A comparative studdy of interface reconstruction methods for
  multi-material {ALE} simulations.
\newblock {\em Journal of Computational Physics}, 229:2432--2452, 2010.

\bibitem{hypre2_14}
Lawrence Livermore~National Lab.
\newblock hypre version 2.14.

\bibitem{lewis_book}
E.~E. Lewis and W.~F. Miller.
\newblock {\em Computational Methods of Neutron Transport}.
\newblock American Nuclear Society, La Grange Park, IL, 1993.

\bibitem{lingus}
C.~Lingus.
\newblock Analytical tests cases for neutrons and radiation transport codes.
\newblock In {\em Second Conference on Transport Theory}, pages 655--659, Los
  Alamos, New Mexico, January 1971. United States Atomic Energy Commision-
  Division of Technical Information.

\bibitem{air1}
T~A Manteuffel, S~M\"unzenmaier, J~W Ruge, and B~S Southworth.
\newblock Nonsymmetric {R}eduction-based {A}lgebraic {M}ultigrid.
\newblock {\em SIAM Journal on Scientific Computing}, (submitted), June 2018.

\bibitem{air2}
T~A Manteuffel, J~W Ruge, and B~S Southworth.
\newblock {Nonsymmetric Algebraic Multigrid Based on Local Approximate Ideal
  Restriction ($\ell$AIR)}.
\newblock {\em SIAM Journal on Scientific Computing}, (accepted), Sep. 2018.

\bibitem{McLendon2001}
William~III McLendon, Bruce Hendrickson, Steve Plimpton, and Lawrence
  Rauchwerger.
\newblock {Finding Strongly Connected Components in Parallel Transport Sweeps}.
\newblock In {\em SPAA 2001}, pages 328--329, 2001.

\bibitem{morel_radtran}
J.~E. Morel, T.~A. Wareing, and K.~Smith.
\newblock A linear-discontinuous spatial differencing scheme for {$S_N$}
  radiative transfer calculations.
\newblock {\em Journal of Computational Physics}, 128:445--462, 1996.

\bibitem{Pautz2002}
Shawn Pautz.
\newblock {An Algorithm for Parallel Sn sweeps on Unstructured Meshes}.
\newblock {\em Nuclear Science and Engineering}, 140:111--136, 2002.

\bibitem{Pautz2017}
Shawn~D Pautz and Teresa~S Bailey.
\newblock {Parallel Deterministic Transport Sweeps of Structured and
  Unstructured Meshes with Overloaded Mesh Decompositions}.
\newblock {\em Nuclear Science and Engineering}, 185(1):70--77, 2017.

\bibitem{Plimpton2005}
Steven~J Plimpton, Bruce Hendrickson, Shawn~P Burns, William McLendon, and
  Lawrence Rauchwerger.
\newblock {Parallel Sn sweeps on unstructured grids: Algorithms for
  prioritization, grid partitioning, and cycle detection}.
\newblock {\em Nuclear Science and Engineering}, 150(3):267--283, 2005.

\bibitem{numerical_book}
W.~H. Press, S.~A. Teukolsky, W.~T. Vetterling, and B.~P. Flannery.
\newblock {\em Numerical Recipes: The Art of Scientific Programming}.
\newblock Cambridge University Press, New York, 3rd edition, 2007.

\bibitem{mms}
K.~Salari and P.~Knupp.
\newblock Code verification by the method of manufactured solutions.
\newblock Technical Report {SAND}2000-1444, Sandia National Labs, 2000.

\bibitem{tarjan}
R.~Tarjan.
\newblock Depth-first search and linear graph algorithms.
\newblock {\em {SIAM} Journal on Computing}, 1(2):146--160, 1972.

\bibitem{sn_on_tets}
T.~A. Wareing, J.~M. McGhee, J.~E. Morel, and S.~D. Pautz.
\newblock Discontinuous finite element $s_n$ methods on three-dimensional
  unstructured grids.
\newblock {\em Nuclear Science and Engineering}, 138:256--268, 2001.

\bibitem{woodsPalmer2}
D.~Woods.
\newblock High-order finite elments $s_n$ transport in $x-y$ geometry on meshes
  with curved surfaces in the thick diffusion limit.
\newblock Master's thesis, Oregon State University, 2016.

\bibitem{woodsPalmer1}
D.~Woods, T.~A. Brunner, and T.~S. Palmer.
\newblock High-order finite elements $s_n$ transport in $x-y$ geometry on
  meshses with curved surfaces.
\newblock {\em Transactions of the American Nuclear Society}, 114:377--380,
  2016.

\bibitem{Woods2017a}
D.~Woods and T.~Palmer.
\newblock {Diffusion Synthetic Acceleration for High Order SN Transport on
  Meshes with Curved Surfaces}.
\newblock In {\em Trans. of the Am. Nucl. Soc. (June 2017 ANS meeting)}, 2017.

\end{thebibliography}

\end{document}